\documentclass[11pt,a4paper]{article}

\usepackage{epic,eepic,epsfig,amsmath,amssymb,amsthm}
\usepackage{t1enc}

\usepackage{color}

\usepackage{bbm}

\def\virgp{\raise 2pt\hbox{,}}

\renewcommand{\geq}{\geqslant}
\renewcommand{\leq}{\leqslant}
\def\N{{\mathbb N}}
\def\C{\mathbb{C}}
\def\R{{\mathbb R}}

\def\virgp{\raise 2pt\hbox{,}}
\def\cdotpv{\raise 2pt\hbox{;}}

\def\1{\mathbbm{1}}
\newcommand{\ds}{\displaystyle}

\newtheorem{theorem}{Theorem}[section]

\newtheorem{proposition}[theorem]{Proposition}

\newtheorem{pte}[theorem]{Property}

\theoremstyle{remark}
\newtheorem{remark}{Remark}[section]

\theoremstyle{definition}
\newtheorem{definition}{Definition}[section]

\newtheorem*{notation}{Notation}

\theoremstyle{definition}

\theoremstyle{definition}

\begin{document}

\title{\textbf{A note on energy forms on fractal domains}}

\author{\LARGE{\textbf{Claire David}}}

\maketitle
\centerline{Sorbonne Universit\'es, UPMC Univ Paris 06}

\centerline{CNRS, UMR 7598, Laboratoire Jacques-Louis Lions, 4, place Jussieu 75005, Paris, France}



\maketitle

\section{Introduction}

  In~\cite{Kigami1989},~\cite{Kigami1993},~J. Kigami has laid the foundations of what is now known as analysis on fractals, by allowing the construction of an operator of the same nature of the Laplacian, defined locally, on graphs having a fractal character. The Sierpi\'nski gasket stands out of the best known example. It has, since then, been taken up,  developed and popularized by R.~S.~Strichartz~\cite{Strichartz1999}, \cite{StrichartzLivre2006}.\\

The Laplacian is obtained through a weak formulation, obtained by means of Dirichlet forms, built by induction on a sequence of graphs that converges towards the considered domain. It is these Dirichlet forms that enable one to obtain energy forms on this domain. \\

Yet, things are not that simple. If, for domains like the Sierpi\'nski gasket, the Laplacian is obtained in a quite natural way, one must bear in mind that Dirichlet forms solely depend on the topology of the domain, and not of its geometry. Which means that, if one aims at building a Laplacian on a fractal domain, the topology of which is the same as, for instance, a line segment, one has to find a way of taking account a very specific geometry. We came across that problem in our work on the graph of the Weierstrass function~\cite{ClaireGB}. The solution was thus to consider energy forms more sophisticated than classical ones, by means of normalization constants that could, not only bear the topology, but, also, the very specific geometry of, from now on, we will call~${\cal W}-$curves.\\

It is interesting to note that such problems do not seem to arise so much in the existing literature. In a very general way, one may refer to~\cite{UtaFreiberg2004}, where the authors build an energy form on non-self similar closed fractal curves, by integrating the Lagrangian on this curve.\\

We presently aim at investigating the links between energy forms and geometry. We have chosen to consider fractal curves, specifically, the Sierpi\'nski arrowhead curve, the limit of which is the Sierpi\'nski gasket. Does one obtain the same Laplacian as for the triangle ? The question appears as worth to be investigated.

\newpage

\section{Framework of the study}

We place ourselves, in the following, in the Euclidean plane of dimension~2, referred to a direct orthonormal frame. The usual Cartesian coordinates are~$(x,y)$.\\

\vskip 1cm

\begin{notation}
Given a point~$X\,\in\,\R^2$, we will denote by:
\begin{enumerate}
\item[i.] ${\cal S}im_{X, \frac{1}{2},\frac{\pi}{3}}$ the similarity of ratio~$\displaystyle \frac{1}{2}$, the center of which is~$X$, and the angle,~$\displaystyle \frac{\pi}{3} $ ;
\item[ii.] ${\cal S}im_{X,\frac{1}{2},- \frac{\pi}{3}}$ the similarity of ratio~$\displaystyle \frac{1}{2}$, , the center of which is~$X$, and the angle,~$-\displaystyle \frac{\pi}{3} $.

\end{enumerate}
\end{notation}

\vskip 1cm

\begin{definition}
\noindent Let us consider the following points of~$\R^2$:

$$A=(0,0) \quad , \quad  D=(1,0) \quad , \quad B={\cal S}im_{X, \frac{1}{2},\frac{\pi}{3}}(D)\quad , \quad C={\cal S}im_{X, \frac{1}{2},-\frac{\pi}{3}} (A)  $$

\noindent We will denote by~$V_1$ the ordered set, of the points:

$$\left \lbrace A,B,C,D\right \rbrace$$

\noindent The set of points~$V_1$, where~$A$ is linked to~$B$,~$B$ is linked to~$C$, and where~$C$ is linked to~$D$, constitutes an oriented graph, that we will denote by~$  {\cal SG}^{C}_1 $.~$V_1$ is called the set of vertices of the graph~$ {\cal SG}^{\cal C}_1$.\\

\noindent Let us build by induction the sequence of points:

$$\left ( V_m \right)_{m\in\N^\star}=\left ( X_j^m \right)_{1\leq j \leq {\cal N}^{\cal S}_m, \, m\in\N^\star} \quad , \quad {\cal N}^{\cal S}_m\,\in\,\N^\star $$

\noindent such that:

$$X_1^1=A \quad , \quad  X_2^1=B \quad , \quad  X_3^1=A \quad , \quad X_4^1=D$$

\noindent and for any integers~$m \geq 2$,~$0 \leq j \leq {\cal N}^{\cal S}_m$,~$k\,\in\,\N$,~$\ell\,\in\,\N$:

$$X_{j+k}^m=X_j^{m-1}   \quad \text{if} \quad k \equiv 0 \, [3] $$

$$X_{j+k+\ell}^m={\cal S}im_{X_{j }^{m-1},\frac{1}{2},(-1)^{m+1+\ell }\,\frac{\pi}{3}}\left (X_{j+1}^{m-1} \right)   \quad \text{if} \quad k \equiv 1 \, [3]\quad  \text{and} \quad \ell \,\in\,2\,\N$$

$$X_{j+k+\ell}^m={\cal S}im_{X_{j+1}^{m-1}, \frac{1}{2}, (-1)^{m+\ell}\,\frac{\pi}{3}}\left (X_{j }^{m-1} \right)   \quad \text{if} \quad k \equiv 2 \, [3] \quad  \text{and} \quad \ell \,\in\,\N\setminus 2\,\N$$

\noindent The set of points~$V_m$, where two consecutive points are linked, is an oriented graph, which we will denote by~$  {\cal SG}^{\cal C}_m $.~$V_m$ is called the set of vertices of the graph~${\cal SG}^{\cal C}_m $.

\end{definition}
\vskip 1cm

\begin{pte}
\noindent For any strictly positive integer~$m$:
$$V_m \subset V_{m+1} $$

\end{pte}

\vskip 1cm

\begin{pte}
\noindent If one denotes by~$ \left ({\cal SG}_m \right)_{m \in\N}$ the sequence of graphs that approximate the Sierpi\'{n}ski gasket~$\cal SG$, then, for any strictly positive integer~$m$:

$${\cal SG}^{\cal C}_m \subsetneq {\cal SG}_m$$

\end{pte}

\vskip 1cm

\begin{definition}\textbf{Sierpi\'{n}ski arrowhead curve}\\
\noindent We will denote by~${\cal SG}^{\cal C}$ the limit:

$${\cal SG}^{\cal C} =\displaystyle \lim_{m \to + \infty} {\cal SG}^{\cal C}_m$$

\noindent which will be called the \textbf{Sierpi\'{n}ski arrowhead curve}.

\end{definition}

\vskip 1cm

\begin{pte}
Let us denote by~$\cal SG$ the Sierpi\'{n}ski Gasket. Then:

$$ \displaystyle \lim_{m \to + \infty} {\cal SG}^{\cal C}_m={\cal SG}^{\cal C} = {\cal SG} $$

\end{pte}

\vskip 1cm

\begin{remark}
The sequence of graphs~$ \left ({\cal SG}_m \right)_{m \in\N^\star}$ can also be seen as a \emph{\textbf{Lindenmayer system}} \mbox{("L-system")}, i.e. a set~$(V,\omega,P)$, where~$V$ denotes an alphabet (or, equivalently, the set of constant elements and rules, and variables),~$\omega$, the initial state (also called "axiom"), and~$P$, the production rules, which are to be applied, iteratively, to the initial state.\\
\noindent In the case of the Sierpi\'{n}ski arrowhead curve, if one denotes by:
\begin{enumerate}
\item[\emph{i}.] $F$ the rule: "Draw forward, on one unit length" ;

\item[\emph{ii}.] $+$ the rule: "Turn left, with an angle of~$\displaystyle \frac{\pi}{3} $" ;
\item[\emph{iii}.] $-$ the rule: "Turn right, with an angle of~$\displaystyle \frac{\pi}{3} $" ;

\end{enumerate}

\noindent then:
\begin{enumerate}
\item[\emph{i}.] the variables can be denoted by~$  X$ and~$ Y $ ;

\item[\emph{ii}.] the constants are~$F$, $+$, $-$ ;
\item[\emph{iii}.] the initial state is~$XF$ ;
\item[\emph{iv}.] the production rules are:

$$X \to YF+XF+ Y \quad , \quad  Y \to XF-YF- X  $$
\end{enumerate}
\end{remark}

 \begin{figure}[h!]
 \center{\psfig{height=8cm,width=10cm,angle=0,file=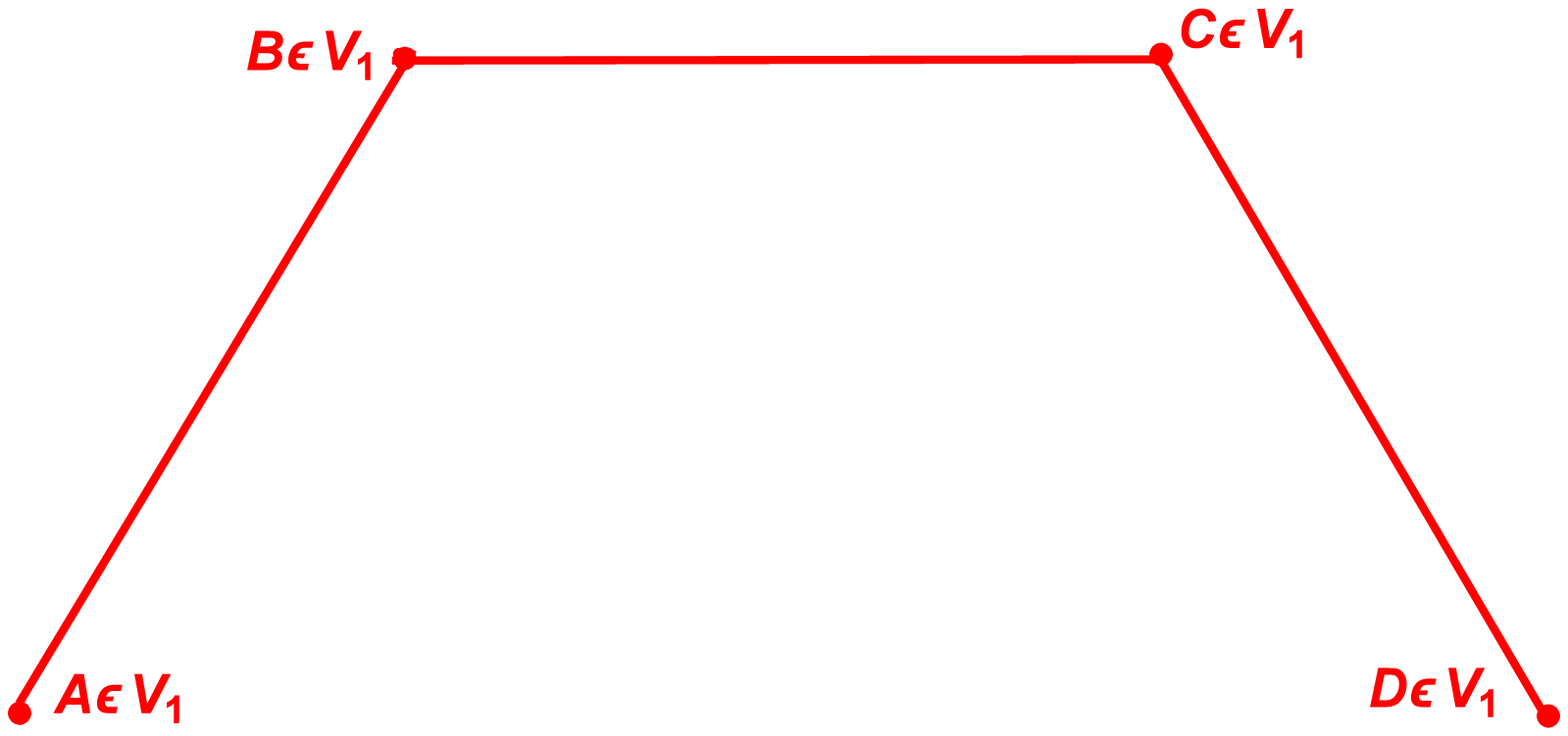}}\\
\caption{The graph~$  {\cal SG}^{\cal C}_1 $.   }
 \end{figure}

 \begin{figure}[h!]
 \center{\psfig{height=8cm,width=10cm,angle=0,file=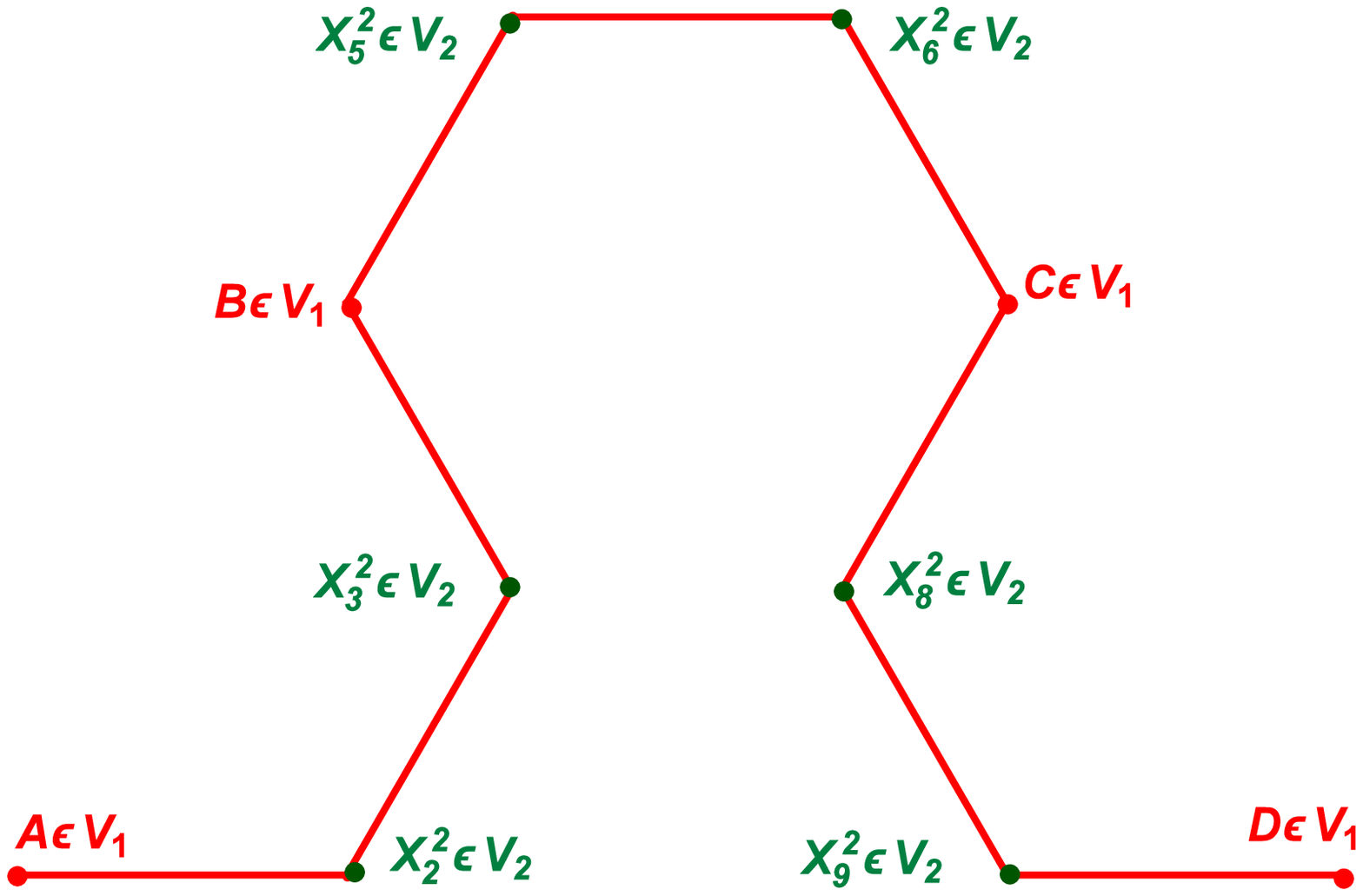}}\\
\caption{The graph~$  {\cal SG}^{\cal C}_2 $.   }
 \end{figure}

 \begin{figure}[h!]
 \center{\psfig{height=8cm,width=10cm,angle=0,file=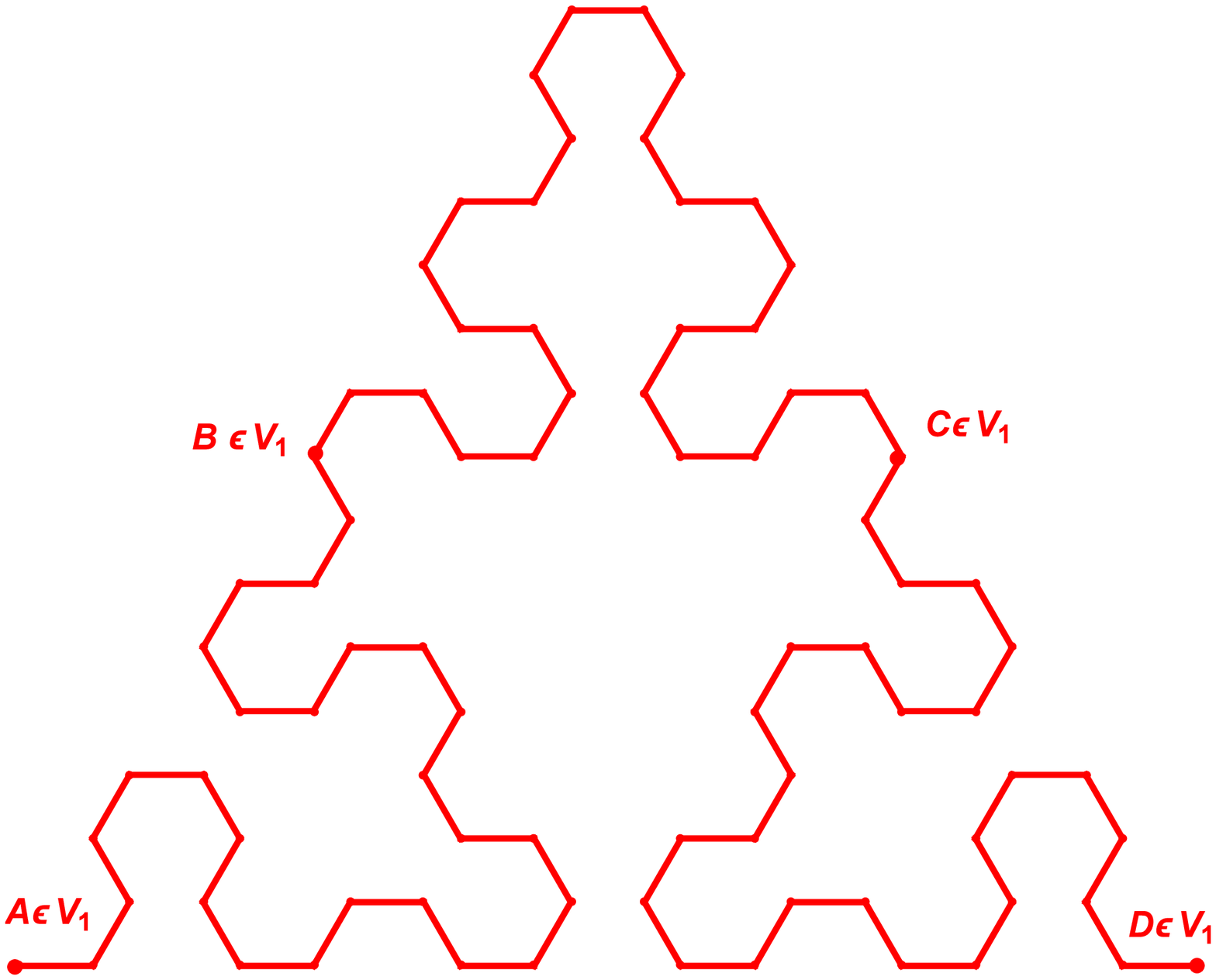}}\\
\caption{The graph~$  {\cal SG}^{\cal C}_4 $.   }
 \end{figure}

 \begin{figure}[h!]
 \center{\psfig{height=8cm,width=10cm,angle=180,file=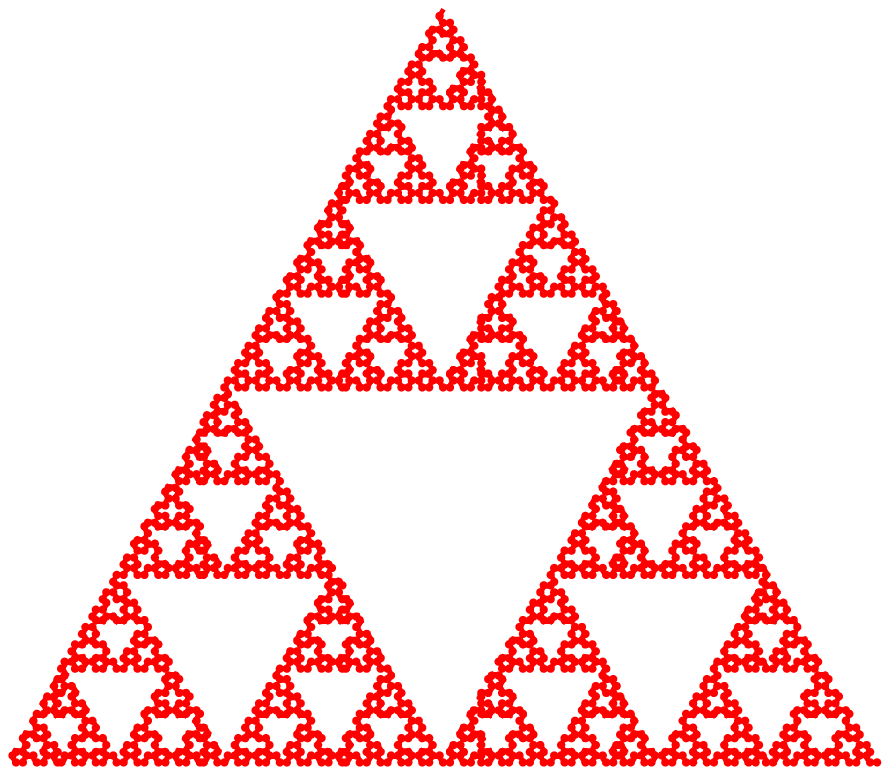}  }\\
\caption{The graph~$  {\cal SG}^{\cal C}_7 $.   }
 \end{figure}

\vskip 1cm

\newpage
\begin{notation}
Given a point~$X\,\in\,\R^2$, we will denote by~${\cal H}_{X, \frac{1}{2}} $ the homothecy of ratio~$\displaystyle \frac{1}{2}$, the center of which is~$X$,.

\end{notation}

\vskip 1cm

\begin{pte}\textbf{Self-similarity properties of the Sierpi\'{n}ski arrowhead curve}\\

\noindent Let us denote by~$E$ the point of~$\R^2$ such that~$A$,~$D$ and~$E$ are the consecutive vertices of a direct equilateral triangle. One may note that~$A$,~$D$ and~$E$ are, also, the frontier vertices of the Sierpi\'{n}ski gasket~$\cal SG$.\\
\noindent The Sierpi\'{n}ski arrowhead curve is self similar with the three homothecies:

$${\cal H}_1={\cal H}_{A, \frac{1}{2}} \quad ,\quad  {\cal H}_2={\cal H}_{D, \frac{1}{2}}  \quad ,\quad {\cal H}_3={\cal H}_{E, \frac{1}{2}} $$

\end{pte}

 \begin{figure}[h!]
 \center{\psfig{height=9cm,width=10cm,angle=0,file=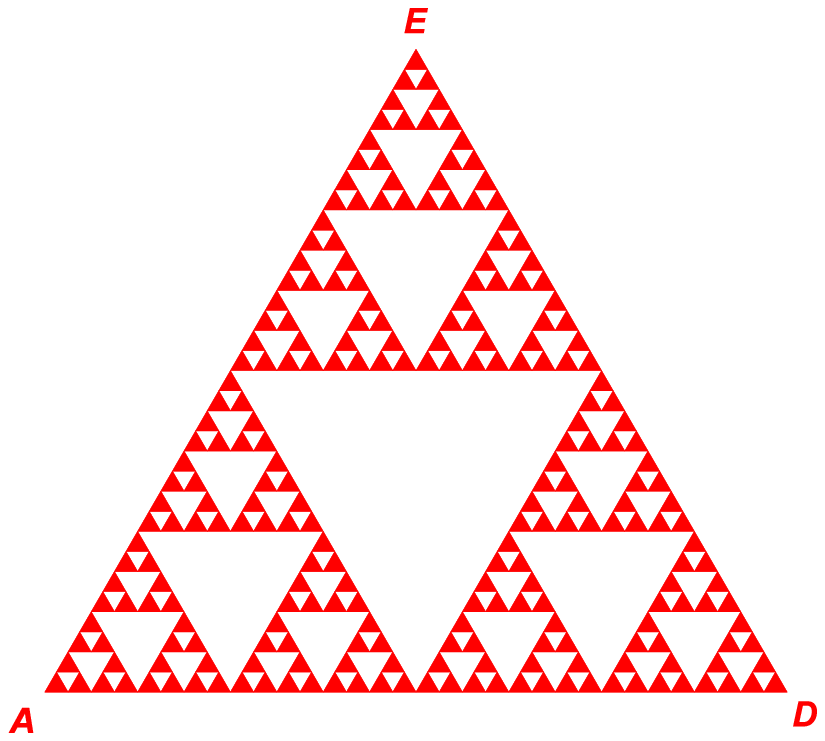}}\\
\caption{The points~$A$,~$D$ and~$E$, as frontier vertices of the Sierpi\'{n}ski gasket.}
 \end{figure}

\vskip 1cm

\begin{proof}The result comes from the self-similarity of the Sierpi\'{n}ski Gasket with respect to those homotecies:

$${\cal SG}=    \underset{  i=1}{\overset{3}{\bigcup}}\,{\cal H}_{i}({\cal SG})$$

\end{proof}
\vskip 1cm

\begin{pte}

\noindent The sequence~$\left ({\cal N}^{\cal S}_m \right)_{m\in\N }$ is an arithmetico-geometric one, with~${\cal N}^{\cal S}_1=4$ as first term:

$$\forall\,m\,\in\,\N \, : \quad {\cal N}^{\cal S}_{m+1}=4\, \left ({\cal N}^{\cal S}_m-1\right)- \left ({\cal N}^{\cal S}_m-1\right)
=3\,{\cal N}^{\cal S}_m-2 $$

\noindent This leads to:

$$\forall\,m\,\in\,\N^\star \, : \quad {\cal N}^{\cal S}_{m+1}=3^m\,\left ( {\cal N}_1^{\cal S}-1 \right)+ 1=3^{m+1}+ 1 $$

\end{pte}
\vskip 1cm

\vskip 1cm

\begin{definition}\textbf{Consecutive vertices on the graph~$ {\cal SG}^{\cal C} $ }\\

\noindent Two points~$X$ and~$Y$ of~${\cal SG}^{\cal C}$ will be called \textbf{\emph{consecutive vertices}} of the graph~${\cal SG}^{\cal C} $ if there exists a natural integer~$m$, and an integer~$j $ of~\mbox{$\left \lbrace  1,...,{\cal N}^{\cal S}_m-1  \right \rbrace$}, such that:

$$X= X_{j }^m   \quad \text{and} \quad Y= X_{j+1 }^m$$

\noindent or:

$$Y= X_{j }^m   \quad \text{and} \quad X= X_{j+1 }^m$$

\end{definition}

\vskip 1cm

\begin{definition}
\noindent For any positive integer~$m$, the~$ {\cal N}^{\cal S}_m$ consecutive vertices of the graph~$  {\cal SG}^{\cal C}_m  $ are, also, the vertices of~$3^{m-1}$
trapezes~${\cal T}_{m,j}$,~\mbox{$1 \leq j \leq 3^{m-1} $}. For any integer~$j$ such that~\mbox{$1 \leq j \leq 3^{m-1} $}, one obtains each trapeze by linking the point number~$j$ to the point number~$j+1$ if~\mbox{$j = i \, \text{mod } 4$},~\mbox{$0 \leq i \leq  2$}, and the point number~$j$ to the point number~$j-3$ if~\mbox{$j =-1 \, \text{mod } 4$}. These trapezes generate a Borel set of~$\R^2$.\\

\noindent In the sequel, we will denote by~${\cal T}_1$ the initial trapeze, the vertices of which are, respectively:

$$A \quad , \quad B \quad , \quad C \quad ,\quad   D$$

\end{definition}
\vskip 1cm

 \begin{figure}[h!]
 \center{\psfig{height=7cm,width=12cm,angle=0,file=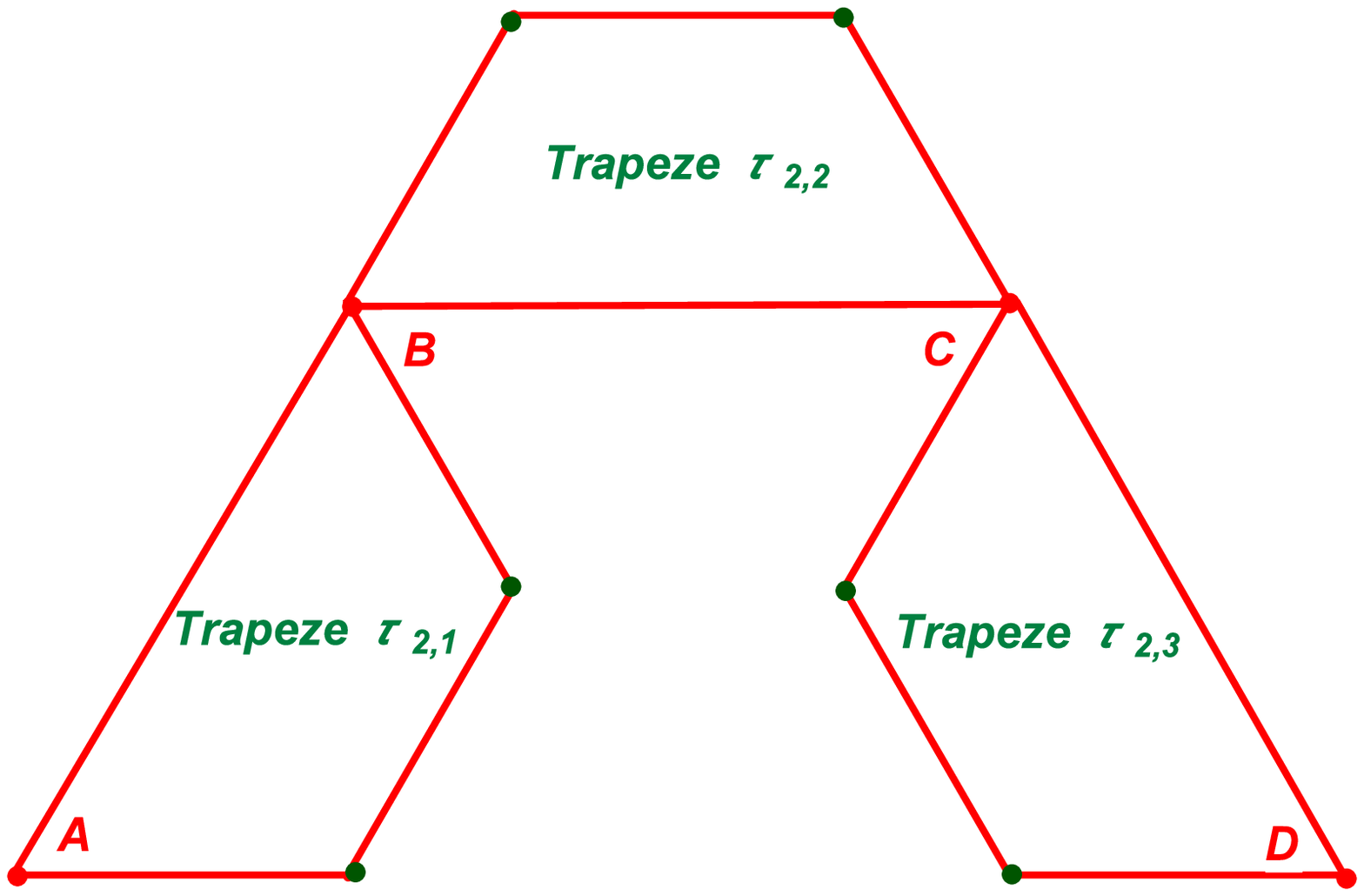}}\\
\caption{The trapezes~${\cal T}_{2,1}$,~${\cal T}_{2,2}$ and~${\cal T}_{2,3}$.}
 \end{figure}

\begin{definition}\textbf{Trapezoidal domain delimited by the graph~$ {\cal SG}^{\cal C}_m  $,~$m\,\in\,\N $}\\

\noindent For any natural integer~$m$, well call \textbf{trapezoidal domain delimited by the graph~$  {\cal SG}^{\cal C}_m  $}, and denote
by~\mbox{$ {\cal D} \left ({\cal SG}^{\cal C}_m \right) $}, the reunion of the~$3^{m-1}$ trapezes~${\cal T}_{m,j}$,~\mbox{$1 \leq j \leq 3^{m-1} $}.\\

\end{definition}
\vskip 1cm

\begin{pte}

\noindent Taking into account that the Lebesgue measure of the first trapeze~${\cal T}_1$ is given by:

$${\cal A}_1={\cal A} \left ({\cal T}_1 \right)=\displaystyle \frac{\sqrt{3}}{4}$$

\noindent one obtains, for any natural integer~$m>1$, the Lebesgue measure of a trapeze~${\cal T}_{m,j}$,~\mbox{$1 \leq j \leq 3^{m-1} $} by noticing that each trapeze is, also, the reunion of three equilateral triangles.\\

\noindent Thus, for any natural integer~$m\geq 2$, the Lebesgue measure of a trapeze~${\cal T}_{m,j}$,~\mbox{$1 \leq j \leq 3^{m-1} $} is given by:

$${\cal A}_m={\cal A} \left ({\cal T}_{m,j} \right)= \displaystyle \frac{3\,{\cal A}_1}{4^m}$$

\end{pte}
\vskip 1cm

\begin{definition}\textbf{Trapezoidal domain delimited by the graph~$  {\cal SG}^{\cal C} $ }\\

\noindent We will call \textbf{trapezoidal domain delimited by the graph~$  {\cal SG}^{\cal C} $}, and denote by~\mbox{$ {\cal D} \left ( {\cal SG}^{\cal C}  \right) $}, the limit:
$$ {\cal D} \left ({\cal SG}^{\cal C}  \right)  = \displaystyle \lim_{m \to + \infty} {\cal D} \left ( {\cal SG}^{\cal C}_m \right) $$

\end{definition}

\vskip 1cm

\begin{notation}

\noindent In the sequel, we will denote by~$d_{\R^2}$ the Euclidean distance on~$\R^2$.
\end{notation}
\vskip 1cm

\begin{definition}\textbf{Edge relation, on the graph~${\cal SG}^{\cal C}$}\\

\noindent Given a natural integer~$m$, two points~$X$ and~$Y$ of~${\cal SG}^{\cal C}_m $ will be called \emph{\textbf{adjacent}} if and only if~$X$ and~$Y$ are two consecutive vertices of~${\cal SG}^{\cal C}_m $. We will write:

$$X \underset{m }{\sim}  Y$$

 \noindent Given two points~$X$ and~$Y$ of the graph~${\cal SG}^{\cal C} $, we will say that~$X$ and~$Y$ are \textbf{\emph{adjacent}} if and only if there exists a natural integer~$m$ such that:
$$X  \underset{m }{\sim}  Y$$
\end{definition}

\vskip 1cm

\begin{pte}\textbf{Euclidean distance of two adjacent vertices of~${\cal SG}^{\cal C}_m $,~$m \,\in\, \N$}\\

\noindent Given a natural integer~$m$, and two points~$X$ and~$Y$ of~${\cal SG}^{\cal C}_m $ such that~$X \underset{m }{\sim}  Y$:

$$d_{\R^2}(X,Y)=\displaystyle \frac{1}{2^m}$$

\end{pte}

\vskip 1cm

\begin{pte}

The set of vertices~$\left (V_m \right)_{m \in\N}$ is dense in~$ { \cal SG}^{\cal C} $.

\end{pte}

\vskip 1cm





\begin{definition}\textbf{Measure, on the domain delimited by the graph~$  { \cal SG}^{\cal C} $ }\\

\noindent We will call \textbf{domain delimited by the graph~$  { \cal SG}^{\cal C}  $}, and denote by~\mbox{$ {\cal D} \left ( { \cal SG}^{\cal C}   \right) $}, the limit:
$$ {\cal D} \left ( {\cal SG}^{\cal C}  \right)  = \displaystyle \lim_{n \to + \infty} {\cal D} \left ( { \cal SG}^{\cal C} _m \right) $$

\noindent which has to be understood in the following way: given a continuous function~$u$ on the graph~${\cal SG}^{\cal C}$, and a measure with full support~$\mu$ on~$\R^2$, then:

$$\displaystyle \int_{ {\cal D} \left ( { \cal SG}^{\cal C} \right)} u\,d\mu  = \displaystyle \lim_{m \to + \infty}
\displaystyle \sum_{j=1}^{3^{m-1} }  \displaystyle \sum_{X \, \text{vertex of }{\cal T}_{m,j} }u\left ( X \right) \,\mu \left (  {\cal T}_{m,j}  \right)$$

\noindent We will say that~$\mu$ is a \textbf{measure, on the domain delimited by the graph~$  { \cal SG}^{\cal C}  $}.
\end{definition}
\vskip 1cm

\begin{definition}\textbf{Dirichlet form} (we refer to the paper \cite{Beurling1985}, or the book \cite{Fukushima1994})\\

 \noindent Given a measured space~$(E, \mu)$, a \emph{\textbf{Dirichlet form}} on~$E$ is a bilinear symmetric form, that we will denote by~${\cal E}$,
 defined on a vectorial subspace~$D$ dense in $L^2(E, \mu) $, such that:\\

\begin{enumerate}

\item For any real-valued function~$u$ defined on~$D$ :  ${\cal E}(u,u) \geq 0$.

\item   $D$, equipped with the inner product which, to any pair~$(u,v)$ of~$D \times D $, associates:

 $$  (u,v)_{\cal{E}}  = (u,v)_{L^2(E,\mu)} + {\cal {E}}(u,v)$$

is a Hilbert space.

\item For any real-valued function~$u$ defined on~$D$, if:
$$ u_\star = \min\, (\max(u, 0) , 1) \,\in \,D$$

\noindent then : ${ \cal{E}}(u_\star,u_\star)\leq { \cal{E}}(u,u)$ (Markov property, or lack of memory property).

\end{enumerate}

\end{definition}

\vskip 1cm

\begin{definition}\textbf{Dirichlet form, on a finite set} (\cite{Kigami1993})\\

 \noindent Let~$V$ denote a finite set~$V$, equipped with the usual inner product which, to any pair~$(u,v)$ of functions defined on~$V$, associates:

  $$(u,v)= \displaystyle \sum_{p\in  V} u(p)\,v(p)$$

  \noindent A \emph{\textbf{Dirichlet form}}on~$V$ is a symmetric bilinear form~${\cal E}$, such that:\\

\begin{enumerate}

\item For any real valued function~$u$ defined on~$V$:  ${\cal E}(u,u) \geq 0$.

\item   $  {\cal {E}}(u,u)= 0$ if and only if~$u$ is constant on~$V$.

\item For any real-valued function~$u$ defined on~$V$, if:
$$ u_\star = \min\, (\max(u, 0) , 1)  $$

\noindent i.e. :

$$\forall \,p \,\in\,V \, : \quad u_\star(p)= \left \lbrace \begin{array}{ccc} 1 & \text{if}& u(p) \geq 1 \\u(p) & \text{si}& 0 <u(p) < 1 \\0  & \text{if}& u(p) \leq 0 \end{array} \right.$$

\noindent then: ${ \cal{E}}(u_\star,u_\star)\leq { \cal{E}}(u,u)$ (Markov property).

\end{enumerate}

\end{definition}

\vskip 1cm

\begin{notation}
\noindent Let us denote by:

$$D_{{\cal SG}^{\cal C}}=D_{{\cal SG} }=\displaystyle \frac{\ln 3}{\ln 2}$$

\noindent the box-dimension (equal to the Hausdorff dimension), of the Sierpi\'{n}ski arrow curve~\mbox{$ {\cal SG}^{\cal C}$}.\\
\noindent For the sake of simplicity, we will from now on denote it by~$D_{{\cal SG} }$.

\end{notation}

\vskip 1cm
Let us now consider the problem of energy forms on our curve. The following points have to be taken into account:\\
\begin{enumerate}
\item[\emph{i}.] As mentioned in the preamble of this work, Dirichlet forms solely depend on the topology of the sequence of graphs that approximate our curve.
\item[\emph{ii}.] Our curve is, indeed, self-similar, yet, it cannot be obtained by means of an iterated function system, as it is the case with the Sierpi\'{n}ski gasket, or the~${\cal W }-$curve we studied in~\cite{ClaireGB}. \\

\end{enumerate}
\noindent Such a problem was studied by~U.~Mosco~\cite{UmbertoMosco2002}, who specifically considered the case of what he called "the Sierpi\'{n}ski curve", or "Sierpi\'{n}ski string". Yet, he did not dealt with the curve itself, but with the Sierpi\'{n}ski gasket: "2D branches (...) meet together". Contrary to the arrow curve, the Sierpi\'{n}ski gasket exhibits self-similarity properties which turn it into a post-critically finite fractal (pcf fractal).\\

Yet, one can find interesting ideas in the work of~U.~Mosco. For instance, he suggests to generalize Riemaniann models to fractals and relate the fractal analogous of gradient forms, i.e. the Dirichlet forms, to a metric that could reflect the fractal properties of the considered structure. The link is to be made by means of specific energy forms.\\

There are two major features that enable one to characterize fractal structures:
\begin{enumerate}
\item[\emph{i}.] Their topology, i.e. their ramification.
\item[\emph{ii}.] Their geometry.\\

\end{enumerate}

The topology can be taken into account by means of classical energy forms (we refer to~\cite{Kigami1989}, ~\cite{Kigami1993},~\cite{Strichartz1999}, \cite{StrichartzLivre2006}).\\
As for the geometry, again, things are not that simple to handle.~U.~Mosco introduces a strictly positive parameter,~$\delta$, which is supposed to reflect the way ramification - or the iterative process that gives birth to the sequence of graphs that approximate the structure - affects the initial geometry of the structure. For instance, if $m$ is a natural integer,~$X$ and~$Y$ two points of the initial graph~$V_1$, and~$\cal M$ a word of length~$m$, the Euclidean distance~$d_{\R^2}(X,Y)$ between ~$X$ and~$Y$ is changed into the effective distance:

$$\left ( d_{\R^2}(X,Y)\right )^\delta  $$

This parameter~$\delta$ appears to be the one that can be obtained when building the effective resistance metric of a fractal structure~(see~\cite{StrichartzLivre2006}), which is obtained by means of energy forms. To avoid turning into circles, this means:
\begin{enumerate}
\item[\emph{i}.] either working, in a first time, with a value~$\delta_0$ equal to one, and, then, adjusting it when building the effective resistance metric ;
\item[\emph{ii}.] using existing results, as done in~\cite{UtaFreiberg2004}.

\end{enumerate}

In the case of the arrow curve, at a step~$m\,\in\,\N^\star$ of the iteration process, the distance between two adjacent points of~${\cal SG}^{\cal C }_m$ is the same as the one between two adjacent points of the graph~${\cal SG}_m$, and take:

$$\delta=\displaystyle \frac{\ln 5}{\ln 4}$$

\vskip 1cm
\newpage
\begin{definition}\textbf{Energy, on the graph~${ \cal SG}^{\cal C} _m  $,~$m \,\in\,\N$, of a pair of functions}\\

 \noindent Let~$m$ be a natural integer, and~$u$ and~$v$ two real valued functions, defined on the set

 $$V_m = \left \lbrace    X_1^m, \hdots,  X_{{\cal N}_m^{\cal S} }^m \right \rbrace $$

 \noindent of the~${\cal N}_m^{\cal S}$ vertices of~${\cal SG}^{\cal C}_m   $.\\

\noindent We introduce \textbf{the energy, on the graph~${ \cal SG}^{\cal C} _m  $, of the pair of functions~$(u,v)$}, as:

$$\begin{array}{ccc}
  {\cal{E}}_{{\cal SG}^{\cal C}_m  }(u,v)
 &= &  \displaystyle \sum_{i=1}^{{\cal N}_m^{\cal S}-1}  \left (\displaystyle \frac{u \left (X_i^m \right)-u \left (X_{i+1}^m \right)}{d_{\R^2}^{\delta}(X,Y)}\right )\,
 \left (\displaystyle \frac{v \left (X_{i }^m \right)-v \left (X_{i+1}^m \right)}{ d_{\R^2}^{\delta}(X,Y)} \right )\\
  &= &  \displaystyle \sum_{i=1}^{{\cal N}_m^{\cal S}-1} 2^{2\,m \,\delta}\,\left ( u \left (X_i^m \right)-u \left (X_{i+1}^m \right) \right )\,
 \left ( v \left (X_{i }^m \right)-v \left (X_{i+1}^m \right)  \right )\\
 \end{array}
$$

 \noindent For the sake of simplicity, we will write it under the form:

$$ {\cal{E}}_{{\cal SG}^{\cal C}_m  }(u,v)= \displaystyle \sum_{X  \underset{m }{\sim}  Y} 4^{ m \,\delta}\,\left (u(X)-u(Y)\right )\,\left(v(X)-v(Y)\right) $$

\end{definition}

\vskip 1cm

\begin{pte}

 \noindent Given a natural integer~$m$, and a real-valued function~$u$, defined on the set of vertices of~${\cal SG}^{\cal C}_m  $, the map, which, to any pair of real-valued, continuous functions~$(u,v)$ defined on the set~$V_m $ of the~${\cal N}_m$ vertices of~${\cal SG}^{\cal C}_m   $, associates:
$$ {\cal{E}}_{{\cal SG}^{\cal C}_m }(u,v)= \displaystyle \sum_{X  \underset{m }{\sim}  Y} 4^{ m \,\delta}\,\left (u(X)-u(Y)\right )\, \left (v(X)-v(Y)\right ) $$

\noindent is a Dirichlet form on~${ \cal SG}^{\cal C} _m  $.\\
\noindent Moreover:

$${\cal{E}}_{{\cal SG}^{\cal C}_m   }(u,u)=0 \Leftrightarrow u  \text{ is constant}$$

\end{pte}

\vskip 1cm

\begin{proposition}\textbf{Harmonic extension of a function, on the graph of Sierpi\'{n}ski arrow curve - Ramification constant}\\

\noindent For any integer~$m>1$, if~$u$ is a real-valued function defined on~$V_{m-1}$, its \textbf{harmonic extension}, denoted by~$ \tilde{u}$, is obtained as the extension of~$u$ to~$V_m$ which minimizes the energy:

$$  {\cal{E}}_{{\cal SG}^{\cal C}_m }(\tilde{u},\tilde{u})= \displaystyle \sum_{X \underset{m }{\sim} Y} 4^{ m \,\delta}\,(\tilde{u}(X)-\tilde{u}(Y))^2 $$

\noindent The link between~$   {\cal{E}}_{{\cal SG}^{\cal C}_m }$ and~$  {\cal{E}}_{{\cal SG}^{\cal C}_{m-1} }$ is obtained through the introduction of two strictly positive constants~$r_m$ and~$r_{m+1}$ such that:

$$   r_{m }\, \displaystyle \sum_{X \underset{m  }{\sim} Y} 4^{ m \,\delta}\, (\tilde{u}(X)-\tilde{u}(Y))^2 =  r_{m-1}\,4^{ m \,\delta}\, \sum_{X \underset{m-1 }{\sim} Y} (u(X)-u(Y))^2$$

\noindent In particular:

$$    r_{2 }\,4^{ 2\, \delta}\, \displaystyle \sum_{X \underset{1  }{\sim} Y} (\tilde{u}(X)-\tilde{u}(Y))^2 =   r_{1}\,4^{  \delta}\,\displaystyle  \sum_{X \underset{1 }{\sim} Y} (u(X)-u(Y))^2$$

\noindent For the sake of simplicity, we will fix the value of the initial constant:~$r_1=1$. One has then:

$$ {\cal{E}}_{ {\cal SG}^{\cal C}_m }(\tilde{u},\tilde{u})= \displaystyle \frac{1}{ r_{1 }}\,  {\cal{E}}_{ {\cal SG}^{\cal C}_1 }(\tilde{u},\tilde{u})$$

\noindent Let us set:

$$r = \displaystyle \frac{1}{r_{1 }} $$

\noindent and:

$$  {\cal{E}}_{m}(u)= r_m\, \sum_{X \underset{m }{\sim} Y}  4^{ m \,\delta}\,(\tilde{u}(X)-\tilde{u}(Y))^2 $$

\noindent Since the determination of the harmonic extension of a function appears to be a local problem, on the graph~$\Gamma_{{\cal W}_{ m-1}}$, which is linked to the
graph~$ {\cal SG}^{\cal C}_m  $ by a similar process as the one that links~$ {\cal SG}^{\cal C}_2$ to~$ {\cal SG}^{\cal C}_1$, one deduces, for any integer~$m>2$:

$$ {\cal{E}}_{{\cal SG}^{\cal C}_m}(\tilde{u},\tilde{u})= \displaystyle \frac{1}{ r_{1 }}\,  {\cal{E}}_{{\cal SG}^{\cal C}_{m-1}}(\tilde{u},\tilde{u})$$

\noindent By induction, one gets:

$$r_m=r_1^m =r^{-m}=3^{-m}$$


\noindent If~$v$ is a real-valued function, defined on~$V_{m-1}$, of harmonic extension~$ \tilde{v}$, we will write:

$$  {\cal{E}}_{m}(u,v)= r^{-m}\,\displaystyle  \sum_{X \underset{m }{\sim} Y}  4^{ m \,\delta}\,(\tilde{u}(X)-\tilde{u}(Y)) \, (\tilde{v}(X)-\tilde{v}(Y)) $$

\noindent The constant~$r^{-1}$, which can be interpreted as a topological one, will be called \textbf{ramification constant}.\\
\noindent For further precision on the construction and existence of harmonic extensions, we refer to~\emph{\cite{Sabot1987}}.
\end{proposition}

\vskip 1cm

\begin{remark}\textbf{Determination of the ramification constant~$r$}\\

\noindent Let us denote by~$u$ a real-valued, continuous function defined on~$V_1$, and by~$\tilde{u}$ its harmonic extension to~$V_2$.\\
\noindent Let us denote by~$a$,~$b$,~$c$ and~$d$ the values of~$u$ on the four consecutive vertices of~$V_1$ (see the following figure):

$$u(A)=a \quad ,\quad u(B)=b \quad ,\quad u(C)=c \quad ,\quad u(D)=d $$

\noindent and by:
\begin{enumerate}
\item[\emph{i}.] ~$e$ and~$f$ the values of~$\tilde{u}$ on the two consecutive vertices~$E$ and~$F$ that are between~$A$ and~$B$:
$$u(E)=e \quad ,\quad u(F)=f $$
\item[\emph{ii}.]~$g$ and~$h$ the values of~$\tilde{u}$ on the two consecutive vertices~$G$ and~$H$ that are between~$B$ and~$C$:
$$u(G)=g \quad ,\quad u(H)=h $$
\item[\emph{iii}.]~$i$ and~$j$ the values of~$\tilde{u}$ on the two consecutive vertices~$I$ and~$J$ that are between~$C$ and~$D$:
$$u(I)=i \quad ,\quad u(J)=j $$
\end{enumerate}

 \begin{figure}[h!]
 \center{\psfig{height=5cm,width=7cm,angle=0,file=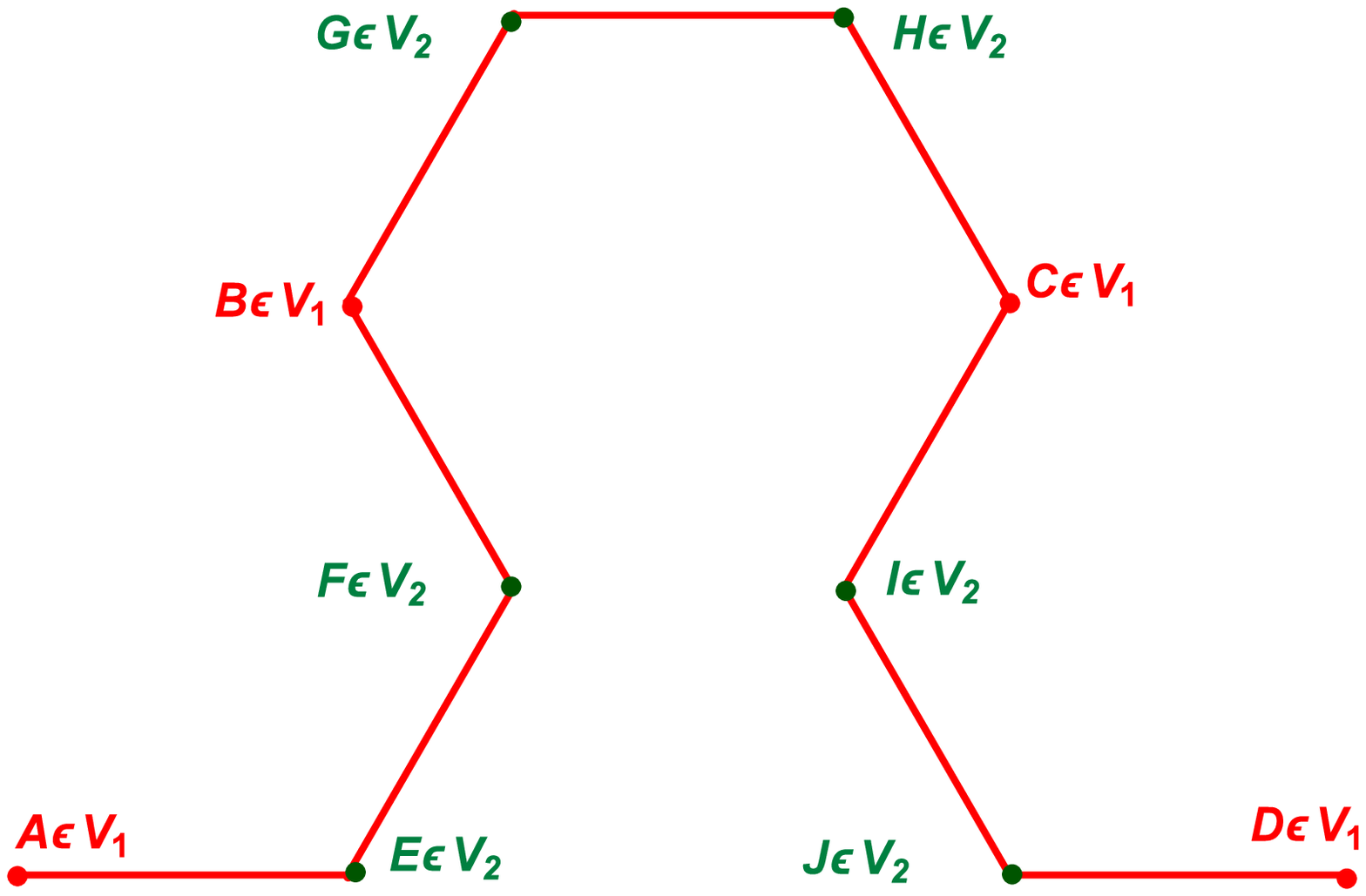}}\\
\caption{Determination of the ramification constant between graphs of level 1 and 2. }
 \end{figure}

\noindent One has:

$$  {\cal{E}}_{{\cal SG}^{\cal C}_1 }(\tilde{u},\tilde{u})= (a-b^2+(b-c)^2+(c-d)^2 $$

$$  {\cal{E}}_{{\cal SG}^{\cal C}_2 }(\tilde{u},\tilde{u})= (a-e)^2+(e-f)^2+(b-f)^2+(g-b)^2+(h-g)^2+(c-h)^2+(i-c)^2+(j-i)^2+(d-j)^2 $$

\noindent Since the harmonic extension~$\tilde{u}$ minimizes~$  {\cal{E}}_{{\cal SG}^{\cal C}_2 }$, the values of~$e$,~$f$,~$g$,~$h$,~$i$,~$j$ are to be found among the critical points~$e$,~$f$,~$g$,~$h$,~$i$,~$j$ such that:
\footnotesize

$$\displaystyle \frac{\partial {\cal{E}}_{{\cal SG}^{\cal C}_2 }(\tilde{u},\tilde{u})}{\partial e}=0  \quad ,\quad
\displaystyle \frac{\partial {\cal{E}}_{{\cal SG}^{\cal C}_2 }(\tilde{u},\tilde{u})}{\partial f}=0  \quad ,\quad
\displaystyle \frac{\partial {\cal{E}}_{{\cal SG}^{\cal C}_2 }(\tilde{u},\tilde{u})}{\partial g}=0  \quad ,\quad
\displaystyle \frac{\partial {\cal{E}}_{{\cal SG}^{\cal C}_2 }(\tilde{u},\tilde{u})}{\partial h}=0  \quad ,\quad
\displaystyle \frac{\partial {\cal{E}}_{{\cal SG}^{\cal C}_2 }(\tilde{u},\tilde{u})}{\partial i}=0  \quad ,\quad
\displaystyle \frac{\partial {\cal{E}}_{{\cal SG}^{\cal C}_2 }(\tilde{u},\tilde{u})}{\partial j}=0  $$

\normalsize

\noindent This leads to:

$$e=\displaystyle \frac{2\, a +  b)}{3} \quad ,\quad f=\displaystyle \frac{2\,(a + 2\,b)}{3}  \quad ,\quad
g=\displaystyle \frac{2\, b +  c }{3}  \quad ,\quad h=\displaystyle \frac{2\,(b + 2\,c)}{3}  \quad ,\quad
i=\displaystyle \frac{2\, c +  d)}{3}  \quad ,\quad j=\displaystyle \frac{2\,(c + 2\,d)}{3} $$

\noindent and:

$$  {\cal{E}}_{{\cal SG}^{\cal C}_2 }(\tilde{u},\tilde{u})=
\displaystyle \frac{1}{3}\,  {\cal{E}}_{{\cal SG}^{\cal C}_2 }(\tilde{u},\tilde{u})$$

\noindent Thus:

$$r^{-1}=\displaystyle \frac{1}{3}$$

\noindent One may note that the ramification constant is exactly equal to one plus the number of points that arise in~$V_{m+1}$, for any value of the strictly positive integer~$m$,  between two consecutive vertices of~$V_m$. We thus fall back on the results we previously obtained in~\cite{DavidRiane},~\cite{ClaireGB} for the graph of the Weierstrass function.

\end{remark}

\vskip 1cm





\begin{definition}\textbf{Energy scaling factor}\\

\noindent By definition, the \textbf{energy scaling factor} is the strictly positive constant~$\rho$ such that, for any integer~$m>1$, and any real-valued function~$u$ defined on~$V_{m }$:

$$  {\cal{E}}_{{\cal SG}^{\cal C}_m }(u,u)= \rho\, {\cal{E}}_{{\cal SG}^{\cal C}_m }\left (u_{\mid V_{m-1}},u_{\mid V_{m-1}} \right)$$

\end{definition}
\vskip 1cm

\begin{proposition}
The energy scaling factor~$\rho$ is linked to the topology and the geometry of the fractal curve by means of the relation:

$$\rho =\displaystyle \frac{4^{\delta}}{3}$$

\end{proposition}

\vskip 1cm







\begin{definition}\textbf{Dirichlet form, for a pair of continuous functions defined on the graph~$  {\cal SG}^{\cal C}$}\\

 \noindent We define the Dirichlet form~$\cal{E}$ which, to any pair of real-valued, continuous functions~$(u,v)$ defined on the Sierpi\'{n}ski arrow curve~$ {\cal SG}^{\cal C}$, associates, subject to its existence:

$$
  {\cal{E}} (u,v)= \displaystyle \lim_{m \to + \infty}  {\cal{E}}_{m }\left (u_{\mid V_m},v_{\mid V_m}\right)=
  \displaystyle \lim_{m \to + \infty}\displaystyle \sum_{X  \underset{m }{\sim}  Y} r^{-m}\,4^{ m \,\delta}\, \left (u_{\mid V_m}(X)-u_{\mid V_m}(Y)\right )\,\left(v_{\mid V_m}(X)-v_{\mid V_m}(Y)\right) $$

\end{definition}

\vskip 1cm

\begin{definition}\textbf{Normalized energy, for a continuous function~$u$, defined on the Sierpi\'{n}ski arrow curve}\\
\noindent Taking into account that the sequence~$\left (\mathcal{E}_m\left ( u_{\mid V_m} \right)\right)_{m\in\N}$ is defined on
$$\ds{V_\star =\underset{{i\in \N}}\bigcup \,V_i}$$

\noindent one defines the \textbf{normalized energy}, for a continuous function~$u$, defined on the curve~$   {\cal SG}^{\cal C  }$, by:

$$\mathcal{E}(u)=\underset{m\rightarrow +\infty}\lim  \mathcal{E}_m \left ( u_{\mid V_m} \right)$$

\end{definition}

\vskip 1cm

\begin{notation}
\noindent We will denote by~$\text{dom}\,{\cal E}$ the subspace of continuous functions defined on~${\cal SG}^{\cal C}$, such that:

$$\mathcal{E}(u)< + \infty$$

\end{notation}

\vskip 1cm

\begin{notation}
\noindent We will denote by~$\text{dom}_1\,{ \cal E}$ the subspace of continuous functions defined on~${\cal SG}^{\cal C}$, which take the value on~$V_1$, such that:

$$\mathcal{E}(u)< + \infty$$

\end{notation}

\vskip 1cm

\newpage

\section{Laplacian of a continuous function, on the Sierpi\'{n}ski arrowhead curve}

\begin{definition}\textbf{Self-similar measure, on the graph of the Sierpi\'{n}ski arrow curve}\\

\noindent A measure~$\mu$ on~$\R^2$ will be said to be \textbf{self-similar} for the domain delimited  by the Sierpi\'{n}ski arrow curve, if there exists a family of strictly positive pounds~\mbox{$\left (\mu_1,\mu_2,\mu_3\right)  $} such that:


$$ \mu= \displaystyle \sum_{i=1}^{3} \mu_i\,\mu\circ {\cal H}_i^{-1} \quad, \quad \displaystyle \sum_{i=1}^{3} \mu_i =1$$

\noindent For further precisions on self-similar measures, we refer to the works of~J.~E.~Hutchinson~(see \cite{Hutchinson1981}).

\end{definition}

\vskip 1cm
\begin{pte}\textbf{Building of a self-similar measure, for the domain delimited by the Sierpi\'{n}ski arrow curve}\\

\noindent The Dirichlet forms mentioned in the above require a positive Radon measure with full support. The choice of a self-similar measure, which is, mots of the time, built with regards to a reference set, of measure~1, appears, first, as very natural. R.~S.~Strichartz~\cite{Strichartz1995},~\cite{Strichartz1999}, showed that one can simply consider auto-replicant measures~$ \tilde{\mu}$, i.e. measures~$ \tilde{\mu}$ such that:

$$ \tilde{\mu}= \displaystyle \sum_{i=1}^{ 3} \tilde{\mu}_i\,\tilde{\mu}\circ {\cal H}_i^{-1} \qquad (\star)$$

\noindent where~$\left (\tilde{\mu}_1, \tilde{\mu}_2, \tilde{\mu}_3\right) $ denotes a family of strictly positive pounds.\\

\noindent This latter approach appears as the best suited in our study, since, in the case of the graph~${\cal SG}^{\cal C}$, the initial set consists of the trapeze~${\cal T}_0$, the measure of which, equal to its surface, is not necessarily equal to~1.\\

\noindent Let us assume that there exists a measure~$\tilde{\mu}$ satisfying~($\star$).\\
\noindent  Relation~$(\star)$ yields, for any set of trapezes~${\cal T}_{m,j}$,~\mbox{$m\,\in\,\N$},~\mbox{$1 \leq j \leq 3^{m-1}$}:

$$ \tilde{\mu} \left ( \underset{1 \leq j \leq 3^{m-1}  }{\bigcup} {\cal T}_{m,j} \right) = \displaystyle \sum_{i=1}^{ 3} \tilde{\mu}_i\, \tilde{\mu}\left ( {\cal H}_i^{-1}   \left ( \underset{1 \leq j \leq 3^{m-1} }{\bigcup}{\cal T}_{m,j} \right) \right)$$

\noindent and, in particular:

$$\tilde{\mu} \left ( {\cal H}_1\left ({\cal T}_1\right)  \cup  {\cal H}_2\left ({\cal T}_1\right) \cup  {\cal H}_3\left ({\cal T}_1\right)  \right)
=\displaystyle \sum_{i=1}^{ 3} \tilde{\mu}_i \,\tilde{\mu}\left (  {\cal T}_1\right)$$

\noindent i.e.:

$$\displaystyle \sum_{i=1}^{ 3}  \tilde{\mu} \left (   {\cal H}_i\left ( {\cal T}_1\right)\right)
=\displaystyle \sum_{i=1}^{ 3} \tilde{\mu}_i \,\tilde{\mu} \left (  {\cal T}_1\right)$$

\noindent The convenient choice, for any~$i$ of~\mbox{$ \left \lbrace 1, 2,3 \right \rbrace$ }, is:
$$\tilde{\mu}_i= \displaystyle \frac{  \tilde{\mu}\left ( {\cal H}_i\left ({\cal T}_1\right) \right) }{\tilde{\mu} \left (  {\cal T}_1\right)}=\displaystyle \frac{3}{4}$$

\noindent One can, from the measure~$\tilde{\mu}$, build the self-similar measure~$\mu $, such that:

$$ \mu= \displaystyle \sum_{i=1}^{3} \mu_i\,\mu\circ {\cal H}_i^{-1}  $$

\noindent where~$\left (\mu_i\right)_{1 \leq i \leq  3}$ is a family of strictly positive pounds, the sum of which is equal to~1.\\

\noindent One has simply to set, for any~$i$ of~\mbox{$ \left \lbrace   1,2,3 \right \rbrace$ }:
$$\mu_i= \displaystyle \frac{ 4\,\tilde{\mu}_i}{9} $$

\noindent The measure~$ \mu $ is self-similar, for the domain delimited by the Sierpi\'{n}ski arrowhead curve.





\end{pte}

\vskip 1cm

\begin{definition}\textbf{Laplacian of order~$m\,\in\,\N^\star$}\\

\noindent For any strictly positive integer~$m$, and any real-valued function~$u$, defined on the set~$V_m$ of the vertices of the graph~$ {SG}^{\cal C}_m $, we introduce the Laplacian of order~$m$,~$\Delta_m(u)$, by:

$$\Delta_m u(X) = \displaystyle\sum_{Y \in V_m,\,Y\underset{m}{\sim} X} \left (u(Y)-u(X)\right)  \quad \forall\, X\,\in\, V_m\setminus V_0 $$

\end{definition}

\vskip 1cm

\begin{definition}\textbf{Harmonic function of order~$m\,\in\,\N^\star$}\\

\noindent Let~$m$ be a strictly positive integer. A real-valued function~$u$,defined on the set~$V_m$ of the vertices of the graph~$ {\cal SG}^{\cal c}_m $, will be said to be \textbf{harmonic of order~$m$} if its Laplacian of order~$m$ is null:

$$\Delta_m u(X) =0 \quad \forall\, X\,\in\, V_m\setminus V_0 $$

\end{definition}

\vskip 1cm

\begin{definition}\textbf{Piecewise harmonic function of order~$m\,\in\,\N^\star$}\\

\noindent  Given a strictly positive integer~$m$, a real valued function~$u$, defined on the set of vertices of~${\cal SG}^{\cal C}$, is said to be \textbf{piecewise harmonic function of order~~$m$} if, for any word~${\cal M}$ of length~$ m$,~$u\circ T_{\cal M}$ is harmonic of order~$m$.

\end{definition}

\vskip 1cm

\begin{definition}\textbf{Existence domain of the Laplacian, for a continuous function on the graph~${\cal SG}^{\cal C}$ } (see \cite{Beurling1985})\\

\label{Lapl}
\noindent We will denote by~$\text{dom}\, \Delta$ the existence domain of the Laplacian, on the graph~${\cal SG}^{\cal C}$, as the set of functions~$u$ of~$\text{dom}\, \mathcal{E}$such that there exists a continuous function on~${\cal SG}^{\cal C}$, denoted~$\Delta \,u$, that we will call \textbf{Laplacian of~$u$}, such that :
$$\mathcal{E}(u,v)=-\displaystyle \int_{{\cal D} \left ( {\cal SG}^{\cal C} \right)} v\, \Delta u   \,d\mu \quad \text{for any } v \,\in \,\text{dom}_1\, \mathcal{E}$$
\end{definition}

\vskip 1cm

\begin{definition}\textbf{Harmonic function}\\

\noindent A function~$u$ belonging to~\mbox{$\text{dom}\,\Delta$} will be said to be \textbf{harmonic} if its Laplacian is equal to zero.
\end{definition}

\vskip 1cm

\begin{notation}

In the following, we will denote by~${\cal H}_0\subset \text{dom}\, \Delta$ the space of harmonic functions, i.e. the space of functions~$u \,\in\,\ \text{dom}\, \Delta$ such that:

$$\Delta\,u=0$$

\noindent Given a natural integer~$m$, we will denote by~${\cal S} \left ({\cal H}_0,V_m \right)$ the space, of dimension~$N_b^m$, of spline functions " of level~$m$", ~$u$, defined on~${\cal SG}^{\cal C}$, continuous, such that, for any word~$\cal M$ of length~$m$,~\mbox{$u \circ T_{\cal M}$} is harmonic, i.e.:

$$\Delta_m \, \left ( u \circ T_{\cal M} \right)=0$$

\end{notation}

\vskip 1cm

\begin{pte}

For any natural integer~$m$:

$${\cal S} \left ({\cal H}_0,V_m \right )\subset  \text{dom }{ \cal E}$$

\end{pte}
\vskip 1cm

\begin{pte}
Let~$m$ be a strictly positive integer,~$X \,\notin\,V_0$ a vertex of the graph~${\cal SG}^{\cal C}$, and~\mbox{$\psi_X^{m}\,\in\,{\cal S} \left ({\cal H}_0,V_m \right)$} a spline  function such that:

$$\psi_X^{m}(Y)=\left \lbrace \begin{array}{ccc}\delta_{XY} & \forall& Y\,\in \,V_m \\
 0 & \forall& Y\,\notin \,V_m \end{array} \right. \quad,  \quad \text{where} \quad    \delta_{XY} =\left \lbrace \begin{array}{ccc}1& \text{if} & X=Y\\ 0& \text{else} &  \end{array} \right.$$

\noindent Then, since~$X\, \notin \,V_0$: $\psi_X^{m}\,\in \,\text{dom}_1\, \mathcal{E}$.



\noindent For any function~$u$ of~$\text{dom}\, \mathcal{E}$, such that its Laplacian exists, definition (\ref{Lapl}) applied to~$\psi_X^{m}$ leads to:

$$\mathcal{E}(u,\psi_X^{m})=\eta_{2-D_{\cal W}   }^{-2}\,\mathcal{E}_m(u,\psi_X^{m})= -r^{-m}\,\eta_{2-D_{\cal W}   }^{-2}\,\Delta_m u(X)=- \displaystyle\int_{{\cal D}\left ({\cal SG}^{\cal C}\right)}  \psi_X^{m}\,\Delta u  \, d\mu \approx -\Delta  u(X)\, \displaystyle\int_{{\cal D} \left ( {\cal SG}^{\cal C} \right)}  \psi_X^{m}\, d\mu$$

\noindent since~$\Delta u$ is continuous on~$ {\cal SG}^{\cal C} $, and the support of the spline function~$\psi_X^{m}$ is close to~$X$:

$$\displaystyle\int_{{\cal D} \left ( {\cal SG}^{\cal C} \right)}  \psi_X^{m}\,\Delta u  \, d\mu \approx -\Delta  u(X)\, \displaystyle\int_{{\cal D} \left ( {\cal SG}^{\cal C} \right)}  \psi_X^{m}\, d\mu$$

\noindent By passing through the limit when the integer~$m$ tends towards infinity, one gets:

$$ \displaystyle \lim_{m \to + \infty} \displaystyle\int_{{\cal D} \left ( {\cal SG}^{\cal C} \right)}  \psi_X^{m}\,\Delta_m u  \, d\mu=
 \Delta  u(X)\,\displaystyle \lim_{m \to + \infty}   \displaystyle\int_{{\cal D} \left ( {\cal SG}^{\cal C} \right)}  \psi_X^{m}\, d\mu$$

\noindent i.e.:

$$  \Delta  u(X)= \displaystyle \lim_{m \to + \infty} r^{-m}\,4^{m\,\delta}\, \left (  \displaystyle\int_{{\cal D} \left ( {\cal SG}^{\cal C} \right)}  \psi_X^{m}\, d\mu  \right)^{-1} \,\Delta_m u(X)\,$$

\end{pte}

\vskip 1cm
\newpage

\section{Explicit determination of the Laplacian of a function~$u$ of~\mbox{$\text{dom}\,\Delta$}}

The explicit determination of the Laplacian of a function~$u$ of~\mbox{$\text{dom}\,\Delta$} requires to know:

$$\displaystyle\int_{{\cal D} \left ( {\cal SG}^{\cal C} \right)}  \psi_X^{m}\, d\mu$$

\noindent As it is explained in~\cite{StrichartzLivre2006}, one has just to reason by analogy with the dimension~1, more particulary, the unit interval~$I=[0,1]$, of extremities~$X_0=(0,0)$, and~$X_1=(1,0)$. The functions~$\psi_{X_1}$ and~$\psi_{X_2}$ such that, for any~$Y$ of~$\R^2$ :

$$\psi_{X_1} (Y)=\delta_{X_1Y} \quad  ,  \quad \psi_{X_2} (Y)=\delta_{X_2Y}   $$

\noindent are, in the most simple way, tent functions. For the standard measure, one gets values that do not depend on~$X_1$, or~$X_2$ (one could, also, choose to fix~$X_1$ and~$X_2$ in the interior of~$I$) :

$$\displaystyle\int_{I}  \psi_{X_1}\, d\mu =\displaystyle\int_{I}  \psi_{X_2}\, d\mu=\displaystyle \frac{1}{2}$$

\noindent (which corresponds to the surfaces of the two tent triangles.) \\

 \begin{figure}[h!]
 \center{\psfig{height=8cm,width=10cm,angle=0,file=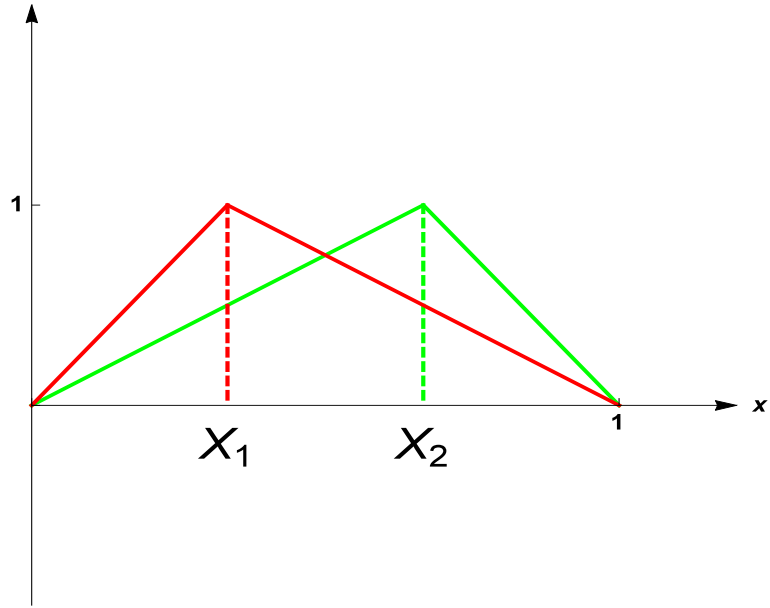}}\\
\caption{The graphs of the spline functions~$\psi_{X_1}$ and~$\psi_{X_2}$.}
 \end{figure}

\noindent In our case, we have to build the pendant, we no longer reason on the unit interval, but on our trapezes. \\
\newpage

\noindent Given a strictly positive integer~$m$, and a vertex~$X$ of the graph~${\cal SG}^{\cal C}_m$, two configurations can occur:\\

\begin{enumerate}

\item[\emph{i}.] the vertex~$X$ belongs to one and only one trapeze~${\cal T}_{m,j}$,~\mbox{$1 \leq j \leq 3^{m-1} $}.\\

In this case, if one considers the spline functions~$\psi_{Z}^{m}$ which correspond to the~$3$ vertices of this trapeze distinct from~$X$:

$$\displaystyle \sum_{Z \, \text{vertex of~${\cal T}_{m,j}$}}  \displaystyle\int_{{\cal D}\left ({\cal SG}^{\cal C}\right)} \psi_{Z}^{m}\,d\mu =\mu\left ({\cal T}_{m,j} \right)$$

\noindent i.e., by symmetry:
$$N_b \,  \displaystyle\int_{{\cal D}\left ({\cal SG}^{\cal C}\right)} \psi_{X}^{m}\,d\mu = \mu\left ({\cal T}_{m,j} \right)$$

\noindent Thus:

$$  \displaystyle\int_{{\cal D}\left ({\cal SG}^{\cal C}\right)}  \psi_{X}^{m}\,d\mu =\displaystyle\frac{1}{4}\, \mu\left ({\cal T}_{m,j} \right)$$

 \begin{figure}[h!]
 \center{\psfig{height=8cm,width=10cm,angle=0,file=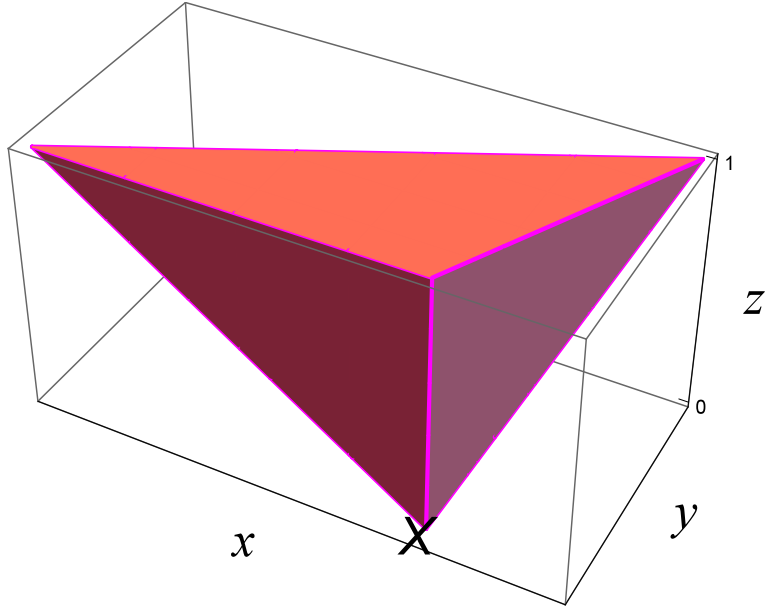}}\\
\caption{The graph of a spline function~$\psi_{X}^{m}$,~$m\,\in\,\N$.}
 \end{figure}

\vskip 1cm

\item[\emph{ii}.] the vertex~$X$ is the intersection point of two trapezes~${\cal T}_{m,j}$ and~${\cal P}_{m,j+1}$,~\mbox{$1 \leq j \leq 3^{m-1}$}.\\

\noindent On has then to take into account the contributions of both trapezes, which leads to:

$$   \displaystyle\int_{{\cal D}\left ({\cal SG}^{\cal C}\right)}  \psi_{X}^{m}\,d\mu =\displaystyle\frac{1}{8}\, \left \lbrace \mu\left ({\cal T}_{m,j} \right)+\mu\left ({\cal T}_{m,j+1} \right) \right\rbrace $$

\end{enumerate}

\newpage

\vskip 1cm
\begin{theorem}

\noindent Let~$u$ be in~\mbox{$\text{dom}\,\Delta$}. Then, the sequence of functions~$\left (f_m \right)_{m\in\N^\star}$ such that, for any strictly positive integer~$m$, and any~$X$ of~\mbox{$V_\star\setminus V_1$} :

 $$f_m(X)=r^{-m}\,4^{m\,\delta}\,\left (\displaystyle \int_{{\cal D} \left ({\cal SG}^{\cal C}\right) }   \psi_{X}^{m}\,d\mu\right)^{-1}\,\Delta_m \,u(X) $$

 \noindent  converges uniformly towards~$\Delta\,u$, and, reciprocally, if the sequence of functions~$\left (f_m \right)_{m\in\N^\star}$ converges uniformly towards a continuous function on~\mbox{$V_\star\setminus V_0$}, then:

 $$u \,\in\, \text{dom}\,\Delta$$
\end{theorem}

\vskip 1cm
\begin{proof}

\noindent Let~$u$ be in~\mbox{$\text{dom}\,\Delta$}. Then:

 $$ r^{-m}\,4^{m\,\delta}\,\left (\displaystyle \int_{{\cal D} \left ({\cal SG}^{\cal C}\right) }  \psi_{X}^{m}\,d\mu\right)^{-1}\,\Delta_m \,u(X)=
\displaystyle \frac{\displaystyle \int_{{\cal D} \left ({\cal SG}^{\cal C}\right) }  \Delta\,u \,\psi_{X}^{m}\,d\mu} {\displaystyle \int_{{\cal D} \left ({\cal SG}^{\cal C}\right) }   \psi_{X}^{m}\,d\mu}   $$

\noindent Since~$u$ belongs to~\mbox{$\text{dom}\,\Delta$}, its Laplacian~$\Delta\,u$ exists, and is continuous on the graph~${\cal SG}^{\cal C}$. The uniform convergence of the sequence~$\left (f_m \right)_{m\in\N}$ follows.\\

\noindent Reciprocally, if the sequence of functions~$\left (f_m \right)_{m\in^\star}$ converges uniformly towards a continuous function on~\mbox{$V_\star\setminus V_1$}, the, for any natural integer~$m$, and any~$v$ belonging to~\mbox{$\text{dom}_1\,{\cal E}$}:

$$\begin{array}{ccc} {\cal{E}}_{m }(u,v)
  &=&  \displaystyle \sum_{(X,Y) \,\in \, V_m^2,\,X  \underset{m }{\sim}  Y} r^{-m}\,4^{m\,\delta}\,\left (u_{\mid V_m}(X)-u_{\mid V_m}(Y)\right )\,\left(v_{\mid V_m}(X)-v_{\mid V_m}(Y)\right) \\
  &=& \displaystyle \sum_{(X,Y) \,\in \, V_m^2,\,X  \underset{m }{\sim}  Y} r^{-m}\,4^{m\,\delta}\,\left (u_{\mid V_m}(Y)-u_{\mid V_m}(X )\right )\,\left(v_{\mid V_m}(Y)-v_{\mid V_m}(X)\right) \\
  &=&- \displaystyle \sum_{X \,\in \,V_m\setminus V_1 } r^{-m}\,4^{m\,\delta}\,\sum_{Y\,\in \,V_m, \, Y  \underset{m }{\sim}  X} v_{\mid V_m}(X)\,\left (u_{\mid V_m}(Y)-u_{\mid V_m}(X)\right )   \\
  & &- \displaystyle \sum_{X \,\in \,  V_1 } r^{-m}\,4^{m\,\delta}\,\sum_{Y\,\in \,V_m ,\, Y  \underset{m }{\sim}  X} v_{\mid V_m}(X)\,\left (u_{\mid V_m}(Y)-u_{\mid V_m}(X)\right )   \\
  &=&- \displaystyle \sum_{X \,\in \,V_m\setminus V_1 } r^{-m}\,4^{m\,\delta}\,v(X)\,\Delta_m \,u(X)  \\
    &=&- \displaystyle \sum_{X \,\in \,V_m\setminus V_1 } v(X)\,   \left (\displaystyle \int_{{\cal D} \left ({\cal SG}^{\cal C}\right) }   \psi_{X}^{m}\,d\mu\right) \,
    r^{-m}\,4^{m\,\delta}\, \left (\displaystyle \int_{{\cal D} \left ({\cal SG}^{\cal C}\right) }   \psi_{X}^{m}\,d\mu\right)^{-1}\, \Delta_m \,u(X)  \\
  \end{array}
  $$

\noindent Let us note that any~$X$ of~$V_m\setminus V_1$ admits exactly two adjacent vertices which belong to~$V_m\setminus V_1$, which accounts for the fact that the sum

 $$\displaystyle \sum_{X \,\in \,V_m\setminus V_1 } r^{-m}\,4^{m\,\delta}\,\sum_{Y\,\in \,V_m\setminus V_1 ,\, Y  \underset{m }{\sim}  X} v(X)\,\left (u_{\mid V_m}(Y)-u_{\mid V_m}(X)\right)$$
 \noindent has the same number of terms as:

 $$ \displaystyle \sum_{(X,Y) \,\in \,(V_m\setminus V_1)^2,\,X  \underset{m }{\sim}  Y} r^{-m}\,4^{m\,\delta}\,\left (u_{\mid V_m}(Y)-u_{\mid V_m}(X)\right )\,\left(v_{\mid V_m}(Y)-v_{\mid V_m}(X)\right)  $$

 \noindent For any natural integer~$m$, we introduce the sequence of functions~$\left (f_m \right)_{m\in\N^\star}$ such that, for any~$X$ of~$V_m\setminus V_1$:

 $$f_m(X)=r^{-m}\,4^{m\,\delta}\,\left (\displaystyle \int_{{\cal D} \left ({\cal SG}^{\cal C}\right) }  \psi_{X}^{m}\,d\mu\right)^{-1}\,\Delta_m \,u(X) $$

 \noindent The sequence~$\left (f_m \right)_{m\in\N^\star}$ converges uniformly towards~$\Delta\,u$. Thus:

$$\begin{array}{ccc} {\cal{E}}_{m }(u,v)
  &=&      -   \displaystyle \int_{{\cal D} \left ({\cal SG}^{\cal C}\right) }  \left \lbrace \displaystyle \sum_{X \,\in \,V_m\setminus V_1 } v_{\mid V_m}(X)\,\Delta\,u_{\mid V_m}(X)\, \psi_{X}^{m}\right \rbrace \,d\mu
  \end{array}
  $$

\end{proof}

\subsection{Spectrum of the Laplacian}

In the following, let~$u$ be in~$\text{dom}\, \Delta$. We will apply the \emph{\textbf{spectral decimation method}} developed by~R.~S.~Strichartz \cite{StrichartzLivre2006}, in the spirit of the works of M.~Fukushima et T.~Shima \cite{Fukushima1994}. In order to determine the eigenvalues of the Laplacian~$\Delta\, u$ built in the above, we concentrate first on the eigenvalues~$\left (-{\Lambda_m}\right)_{m\in\N}$ of the sequence of graph Laplacians~$\left (\Delta_m \,u\right)_{m\in\N}$, built on the discrete sequence of graphs~$\left (\Gamma_{{ \cal W}_m}\right)_{m\in\N}$. For any natural integer~$m$, the restrictions of the eigenfunctions of the continuous Laplacian~$\Delta\,u$ to the graph~$\Gamma_{{ \cal W}_m}$ are, also, eigenfunctions of the Laplacian~$\Delta_m$, which leads to recurrence relations between the eigenvalues of order~$m$ and~$m+1$.

\vskip 1cm

We thus aim at determining the solutions of the eigenvalue equation:

$$-\Delta\,u=\Lambda\,u \quad  \text{on } { \cal SG}^{\cal C }$$

\noindent as limits, when the integer~$m$ tends towards infinity, of the solutions of:

$$-\Delta_m\,u=\Lambda_m\,u \quad  \text{on }V_m\setminus V_0$$

\noindent Let~$m \geq 2$. We consider an eigenfunction~$u_{m-1}$ on~\mbox{$V_{m-1}\setminus V_1$}, for the eigenvalue~$\Lambda_{m-1}$. The aim is to extend~$u_{m-1}$ on~\mbox{$V_m\setminus V_1$} in a function~$u_m$, which will itself be an eigenfunction of~$\Delta_m$, for the eigenvalue~$\Lambda_m$, and, thus, to obtain a recurrence relation between the eigenvalues~$\Lambda_m$ and~$\Lambda_{m-1}$. Given three consecutive vertices of~$ { \cal SG}^{\cal C}_{m-1} $,~$X_k$,~$X_{k+1}$,~$X_{k+2}$, where~$k$ denotes a generic natural integer, we will denote by~$Y_{k+1}$,~$Y_{k+2}$ the points of~\mbox{$V_m\setminus V_{m-1}$} such that:~$Y_{k+1}$,~$Y_{k+2}$ are  between~$X_k$ and~$X_{k+1}$, by~$Y_{k+4}$,~$Y_{k+5}$, the points of~\mbox{$V_m\setminus V_{m-1}$} such that:~$Y_{k+4}$,~$Y_{k+5}$ are between~$X_{k+1}$ and~$X_{k+2}$, and by~$Y_{k+7}$,~$Y_{k+8}$, the points of~\mbox{$V_m\setminus V_{m-1}$} such that:~$Y_{k+7}$,~$Y_{k+8}$ are between~$X_{k+2}$ and~$X_{k+3}$. For the sake of consistency, let us set:

$$Y_{k   }=X_{k } \quad , \quad Y_{k+3  }=X_{k+1} \quad , \quad  Y_{k+6 }=X_{k+2}\quad , \quad  Y_{k+9 }=X_{k+3}$$

 \begin{figure}[h!]
 \center{\psfig{height=5cm,width=8cm,angle=0,file=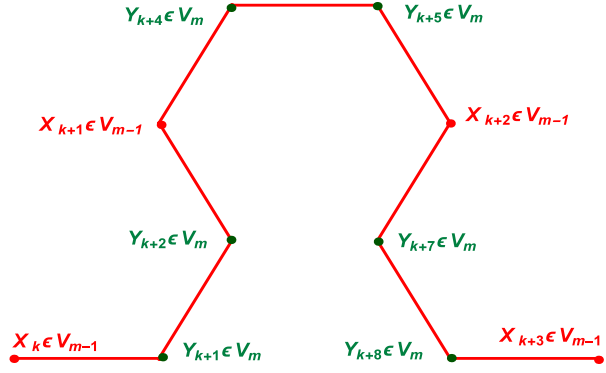}}\\
\caption{The points~$X_k$,~$X_{k+1}$,~$X_{k+2}$,~$X_{k+3}$, and~$Y_{k }$,~$\hdots$,~$Y_{k+9 }$.}
 \end{figure}

\vskip 1cm

\noindent The eigenvalue equation in~$\Lambda_m$ leads to the following system:

$$\left \lbrace
\begin{array}{ccccc}
 \left \lbrace \Lambda_{m }-2 \right\rbrace \,u_m\left ( Y_{k+i+1} \right) &=& - u_{m }\left ( Y_{k+i} \right) -u_m\left ( Y_{k+i+2} \right) &=& - u_{m-1}\left ( X_{k+i} \right) -u_m\left ( Y_{k+i+2} \right) \\
 \left \lbrace \Lambda_{m }-2 \right\rbrace \,u_m\left ( Y_{k+i+2} \right) &=& - u_{m  }\left ( Y_{k+i+1} \right) -u_{m }\left ( Y_{k +i+1} \right)
  &=& - u_{m -1}\left ( X_{k+i+1} \right) -u_{m }\left ( Y_{k +i+1} \right)
\end{array}
\right. \, , \, 0 \leq i \leq 2$$

\noindent The sequence~$\left ( u_m\left ( Y_{k+i} \right) \right)_{0 \leq i \leq 9 }$ satisfies a second order recurrence relation, the characteristic equation of which is:


 $$ r^2+\left \lbrace \Lambda_{m }-2 \right\rbrace \,r+1=0 $$

 \noindent The discriminant is:

$$ \delta_m= \left \lbrace \Lambda_{m }-2 \right\rbrace ^2-4= \omega_m^2 \quad, \quad \omega_m\,\in\,\C$$

\noindent The roots~$r_{1,m}$ and~$r_{2,m}$ of the characteristic equation are the scalar given by:

$$r_{1,m}=\displaystyle \frac{2-\Lambda_{m }-\omega_m}{2}  \quad, \quad r_{2,m}=\displaystyle \frac{2-\Lambda_{m }+\omega_m}{2}$$

\noindent One has then, for any natural integer~$i$ of~\mbox{$\left \lbrace 0,\hdots,9 \right \rbrace $} :

$$ u_m\left ( Y_{k+i } \right) = \alpha_m\, r_{1,m}^i +\beta_m\, r_{2,m}^i$$

\noindent where~$\alpha_m$ and~$\beta_m$ denote scalar constants.

\noindent The extension~$u_m$ of~$u_{m-1}$ to~\mbox{$V_m\setminus V_1$} has to be an eigenfunction of~$\Delta_m$, for the eigenvalue~$\Lambda_m$.\\
\noindent Since~$u_{m-1}$ is an eigenfunction of~$\Delta_{m-1}$, for the eigenvalue~$\Lambda_{m-1}$, the sequence~\mbox{$\left ( u_{m-1}\left ( X_{k+i} \right) \right)_{0 \leq i \leq  9}$} must itself satisfy a second order linear recurrence relation which be the pendant, at order~$m$, of the one satisfied by the sequence~$\left ( u_m\left ( Y_{k+i} \right) \right)_{0 \leq i \leq 9}$, the characteristic equation of which is:

 $$ \left \lbrace \Lambda_{m-1 }-2 \right\rbrace \,r= - 1 -r^2$$

\noindent and discriminant:

$$ \delta_{m-1}= \left \lbrace \Lambda_{m-1 }-2 \right\rbrace ^2-4= \omega_{m-1}^2 \quad, \quad \omega_{m-1}\,\in\,\C$$

\noindent The roots~$r_{1,m-1}$ and~$r_{2,m-1}$ of this characteristic equation are the scalar given by:

$$r_{1,m-1}=\displaystyle \frac{2-\Lambda_{m-1}-\omega_{m-1}}{2} \quad, \quad r_{2,m-1}=\displaystyle \frac{2-\Lambda_{m-1 }+\omega_{m-1}}{2}$$

\noindent For any integer~$i$ of~\mbox{$\left \lbrace 0,\hdots, 9 \right \rbrace $}:

$$ u_{m-1}\left ( Y_{k+i } \right) = \alpha_{m-1}\, r_{1,m-1}^i +\beta_{m-1}\, r_{2,m-1}^i$$

\noindent where~$\alpha_{m-1}$ and~$\beta_{m-1}$ denote scalar constants.

\noindent From this point, the compatibility conditions, imposed by spectral decimation, have to be satisfied:

$$  \left \lbrace \begin{array}{ccc}
u_{m }\left ( Y_{k  } \right)&=&u_{m-1}\left ( X_{k  } \right)  \\
u_{m }\left ( Y_{k+3  } \right)&=&u_{m-1}\left ( X_{k+1 } \right)  \\
u_{m }\left ( Y_{k+6 } \right)&=&u_{m-1}\left ( X_{k+2 } \right)  \\
u_{m }\left ( Y_{k+9 } \right)&=&u_{m-1}\left ( X_{k+3 } \right)  \\
\end{array}\right.$$

\noindent i.e.:

$$  \left \lbrace \begin{array}{ccccc}
 \alpha_{m }  +\beta_{m } &=& \alpha_{m-1}  +\beta_{m-1}  & {\cal C}_{ m}\\
 \alpha_{m }\, r_{1,m }^{3} +\beta_{m }\, r_{2,m }^{3}&=& \alpha_{m-1}\, r_{1,m-1}  +\beta_{m-1}\, r_{2,m-1} & {\cal C}_{1, m} \\
 \alpha_{m }\, r_{1,m }^{6} +\beta_{m }\, r_{2,m }^{6}&=& \alpha_{m-1}\, r_{1,m-1}^2 +\beta_{m-1}\, r_{2,m-1}^2 & {\cal C}_{2, m}\\
 \alpha_{m }\, r_{1,m }^{9} +\beta_{m }\, r_{2,m }^{9}&=& \alpha_{m-1}\, r_{1,m-1}^3 +\beta_{m-1}\, r_{2,m-1}^3 & {\cal C}_{3, m}\\
\end{array}\right. $$

  \noindent where, for any natural integer~$m$,~$\alpha_{m }$ and~$\beta_{m }$ are scalar constants (real or complex).\\

\noindent Since the graph~$ {\cal SG}^{\cal C}_{ m-1} $ is linked to the graph~${\cal SG}^{\cal C}_{ m } $ by a similar process to the one that links~${\cal SG}^{\cal C}_{  2} $ to~${\cal SG}^{\cal C}_{ 1}$, one can legitimately consider that the constants~$\alpha_m$ and~$\beta_m$ do not depend on the integer~$m$:

$$\forall\,m\,\in\,\N^\star \, : \quad \alpha_m=\alpha \,\in\,\R  \quad, \quad \beta_m=\beta \,\in\,\R $$

\noindent The above system writes:

$$  \left \lbrace \begin{array}{ccc}
 \alpha \, r_{1,m }^{3} +\beta \, r_{2,m }^{3}&=& \alpha \, r_{1,m-1}  +\beta \, r_{2,m-1}  \\
 \alpha \, r_{1,m }^{6} +\beta \, r_{2,m }^{6}&=& \alpha \, r_{1,m-1}^2 +\beta \, r_{2,m-1}^2 \\
 \alpha \, r_{1,m }^{8} +\beta \, r_{2,m }^{8}&=& \alpha \, r_{1,m-1}^4 +\beta \, r_{2,m-1}^4 \\
\end{array}\right.$$

\noindent One has then to consider the following configurations:
\begin{enumerate}
\item[\emph{i}.] \underline{First case:}  \\

For any natural integer~$m$ :

$$r_{1,m }\,\in\,\R \quad , \quad r_{2,m }\,\in\,\R$$

\noindent and, more precisely:

$$r_{1,m }<0 \quad , \quad r_{2,m }<0$$

\noindent since the function~$\varphi$, which, to any real number~$x \geq 4$, associates:

$$\varphi(x)=  \displaystyle \frac{2-x+\varepsilon \, \sqrt{\left \lbrace x-2 \right\rbrace ^2-4}}{2}
\quad , \quad \varepsilon \,\in \, \left \lbrace -1, 1 \right \rbrace $$

\noindent is strictly increasing on~$]4,+ \infty[$. Due to its continuity, is is a bijection of~$]4,+ \infty[$ on~\mbox{$ \varphi\left (  ]4,+ \infty[  \right)=]  -1,0 [ $}.\\

\noindent Let us introduce the function~$\phi$, which, to any real number~~$x \geq 2$, associates:

$$\phi(x)= | \varphi(x)|= \displaystyle \frac{-2+x-\varepsilon \,\sqrt{\left \lbrace x-2 \right\rbrace ^2-4}}{2}$$

\noindent where~\mbox{$\varepsilon\,\in\,\left \lbrace -1, 1 \right \rbrace$}.

\noindent The function~$\phi$ is a bijection of~$]4,+ \infty[$ on~\mbox{$ \phi\left (  ]4,+ \infty[  \right)=] 0,1 [ $}.
We will denote by~$\phi^{-1}$ its inverse bijection:

$$\forall\, \,x \,\in\,]0,1[ \, : \quad \phi^{-1}(x)= \displaystyle \frac{(y+1)^2}{y}$$. \\

\noindent One has then:

$$\varphi\left ( \Lambda_{m-1} \right) =\displaystyle \frac{2-\Lambda_{m-1}+\varepsilon\,\omega_{m-1}}{2}  \leq 0  $$

\noindent This yields:

$$
(-1)^{3} \, \left (\varphi\left (  \Lambda_{m } \right) \right)^{3} =   \varphi\left ( \Lambda_{m-1} \right) \leq 0  $$

\noindent which leads to:

$$   \phi\left (  \Lambda_{m } \right) = \left ( \phi\left ( \Lambda_{m-1} \right) \right)^{\frac{1}{3}}  $$






\noindent and:

$$ \Lambda_{m } = \phi^{-1} \left (\left (\phi\left (\Lambda_{m-1 } \right ) \right )^{\frac{1}{3}} \right)
=\displaystyle \frac{\left \lbrace \left (\phi\left (\Lambda_{m-1 } \right ) \right )^{\frac{1}{3}}  +1 \right \rbrace^2}{\left (\phi\left (\Lambda_{m-1 } \right )
\right )^{\frac{1}{3}}   }
= \displaystyle \frac{\left \lbrace \left (\displaystyle \frac{-2+\Lambda_{m-1 }-\varepsilon \,\sqrt{\left \lbrace \Lambda_{m-1 }-2 \right\rbrace ^2-4}}{2}\right)
^{\frac{1}{3}}   +1 \right \rbrace^2}{\left (\displaystyle \frac{-2+\Lambda_{m-1 }-\varepsilon \,\sqrt{\left \lbrace \Lambda_{m-1 }-2 \right\rbrace ^2-4}}{2}\right)
^{\frac{1}{3}}   }
$$

\item[\emph{ii}.] \underline{Second case :}  \\

For any natural integer~$m$:

$$r_{1,m }\,\in\,\C\setminus \R \quad r_{2,m }= \overline{r_{1,m }}\,\in\,\C \setminus \R$$

\noindent Let us introduce:

$$\rho_{ m }= \left | r_{1,m }\right| \,\in\,\R^+ \quad , \quad \theta_m=\text{arg}\, r_{1,m } \quad \text{if} \quad r_{1,m } \neq 0$$

\noindent The above system writes:

$$  \left \lbrace \begin{array}{ccc}
\rho_{ m }^{ 3}\, \left \lbrace \gamma \,\cos \left (3\,\theta_m\right) +\delta\,  \sin \left ( 3\,\theta_m\right) \right \rbrace &=& \rho_{ m-1 } \, \left \lbrace \gamma \,\cos \left (  \theta_{m-1}\right) +\delta\,  \sin \left (  \theta_{m-1}\right) \right \rbrace \\
\rho_{ m }^{6}\, \left \lbrace \gamma \,\cos \left ( 6\,\theta_m\right) +\delta\,  \sin \left (6\,\theta_m\right) \right \rbrace &=&\rho_{ m-1 }^2 \, \left \lbrace \gamma \,\cos \left (  2\,\theta_{m-1}\right) +\delta\,  \sin \left ( 2\,\theta_{m-1}\right) \right \rbrace \\
\rho_{ m }^{9}\, \left \lbrace \gamma \,\cos \left ( 9\,\theta_m\right) +\delta\,  \sin \left (9\,\theta_m\right) \right \rbrace &=&\rho_{ m-1 }^3 \, \left \lbrace \gamma \,\cos \left (  3\,\theta_{m-1}\right) +\delta\,  \sin \left ( 3\,\theta_{m-1}\right) \right \rbrace \\
\end{array}\right.$$

  \noindent where~$\gamma$ and~$\delta$ denote real constants.\\

  \noindent The system is satisfied if:

$$  \left \lbrace \begin{array}{ccc}
\rho_{ m }^{ 3}  &=& \rho_{ m-1 }  \\
 \theta_m  &=&\displaystyle \frac{\theta_{m-1}}{3}
\end{array}\right.$$

\noindent and thus:

$$   \phi\left (  \Lambda_{m } \right) = \left ( \phi\left ( \Lambda_{m-1} \right) \right)^{\frac{1}{N_b}}  $$






\noindent which leads to the same relation as in the previous case:

$$ \Lambda_{m } = \phi^{-1} \left (\left (\phi\left (\Lambda_{m-1 } \right ) \right )^{\frac{1}{3}} \right)
=\displaystyle \frac{\left \lbrace \left (\phi\left (\Lambda_{m-1 } \right ) \right )^{\frac{1}{3}}  +1 \right \rbrace^2}{\left (\phi\left (\Lambda_{m-1 } \right )
\right )^{\frac{1}{3}}   }
= \displaystyle \frac{\left \lbrace \left (\displaystyle \frac{-2+\Lambda_{m-1 }-\varepsilon \,\sqrt{\left \lbrace \Lambda_{m-1 }-2 \right\rbrace ^2-4}}{2}\right)
^{\frac{1}{3}}   +1 \right \rbrace^2}{\left (\displaystyle \frac{-2+\Lambda_{m-1 }-\varepsilon \,\sqrt{\left \lbrace \Lambda_{m-1 }-2 \right\rbrace ^2-4}}{2}\right)
^{\frac{1}{3}}   }
$$

\noindent where~\mbox{$\varepsilon\,\in\,\left \lbrace -1, 1 \right \rbrace$}.

\end{enumerate}

\vskip 1cm

\section{Detailed study of the spectrum of the Laplacian}

 \noindent As exposed by R.~S.~Strichartz in~\cite{StrichartzLivre2006}, one may bear in mind that the eigenvalues can be grouped into two categories:

 \begin{enumerate}

 \item[\emph{i}.]  initial eigenvalues, which a priori belong to the set of forbidden values (as for instance~\mbox{$\Lambda=2$}) ;

 \item[\emph{ii}.] continued eigenvalues, obtained by means of spectral decimation.

 \end{enumerate}

 \vskip 1cm

 \noindent We present, in the sequel, a detailed study of the spectrum of~$\Delta$.
 \vskip 1cm

 \subsection{Eigenvalues and eigenvectors of~$\Delta_2$}

\noindent Let us recall that the vertices of the graph~$ {\cal SG}^{\cal C}_2 $ are:

$$    X_j^2 \quad , \quad  1\leq j \leq 10 $$

\noindent with:

$$    X_1^2 =A \quad , \quad  X_4^2 =B \quad , \quad   X_7^2 =C \quad , \quad  X_{10}^2 =A $$

\noindent For the sake of simplicity, we will set here:

$$    X_2^2 =E \quad , \quad  X_3^2 =F \quad , \quad   X_5^2 =G \quad , \quad  X_6^2 =H \quad , \quad   X_8^2 =I \quad , \quad  X_9^2 =J $$

 \begin{figure}[h!]
 \center{\psfig{height=5cm,width=7cm,angle=0,file=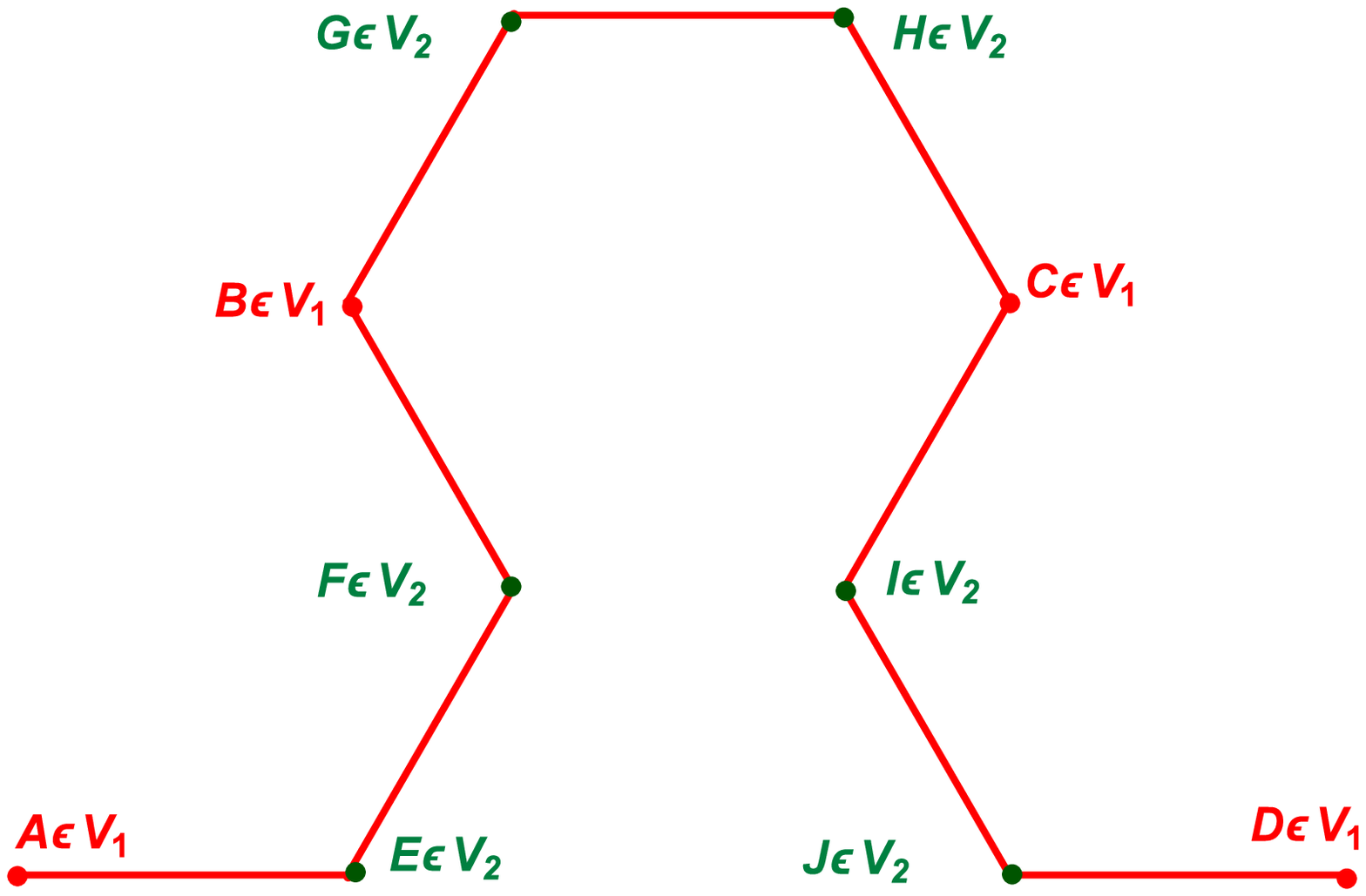}}\\
\caption{Successive values of an eigenfunction on~$V_2$.}
 \end{figure}

\noindent One may note that:

$$\text{Card}\, \left ( V_2\setminus V_1 \right)=10-4=6$$

\noindent Let us denote by~$u$ an eigenfunction, for the eigenvalue~$-\Lambda$. Let us set:

$$u(A)=a\,\in\,\R \quad, \quad u(B)=b\,\in\,\R \quad, \quad  u(C)=c\,\in\,\R \quad, \quad  u(D)=d\,\in\,\R $$

\footnotesize
$$u(E)=e\,\in\,\R \quad, \quad u(F)=f\,\in\,\R \quad, \quad  u(G)=g\,\in\,\R \quad, \quad  u(H)=h\,\in\,\R  \quad, \quad  u(I)=i\,\in\,\R \quad, \quad  u(J)=j\,\in\,\R$$

\normalsize

\noindent One has then:


$$\left \lbrace \begin{array}{ccc}
a +f  &=& -(\Lambda-2)\,e \\
b +e &=& -(\Lambda-2)\,f \\
b +h  &=& -(\Lambda-2)\,g \\
g +c  &=& -(\Lambda-2)\,h \\
c +j &=& -(\Lambda-2)\,i \\
i + d  &=& -(\Lambda-2)\,j \\
\end{array} \right.$$

\noindent One may note that the only "Dirichlet eigenvalues", i.e. the ones related to the Dirichlet problem:

$$u_{|V_1}=0 \quad \text{i.e.} \quad u(A)=u(B)=u(C)=u(D)=0$$

\noindent are obtained for:

$$\left \lbrace \begin{array}{ccc}
 f  &=& -(\Lambda-2)\,e \\
 e &=& -(\Lambda-2)\,f \\
 h  &=& -(\Lambda-2)\,g \\
g    &=& -(\Lambda-2)\,h \\
 j &=& -(\Lambda-2)\,i \\
i    &=& -(\Lambda-2)\,j \\
\end{array} \right.$$

\noindent i.e.:

$$\left \lbrace \begin{array}{ccc}
 f  &=&  (\Lambda-2)^2\,f \\
 e &=&  (\Lambda-2)^2\,e \\
 h  &=&  (\Lambda-2)\,h \\
g    &=&  (\Lambda-2)^2\,g \\
 j &=&  (\Lambda-2)^2\,i \\
i    &=&  (\Lambda-2)^2\,i \\
\end{array} \right.$$

\noindent The forbidden eigenvalue~$\Lambda=2$ cannot thus be a Dirichlet one.\\

\vskip 1cm
\noindent Let us consider the case where:

$$ (\Lambda-2)^2=1$$

\noindent i.e.

$$  \Lambda =1 \quad \text{or} \quad \Lambda =3   $$

\noindent which yields a three-dimensional eigenspace. The multiplicity of the eigenvalue~$\Lambda =1$ is~3.\\\\

\noindent In the same way, the eigenvalue~$\Lambda =3$ yields a three-dimensional eigenspace. the multiplicity of the eigenvalue~$\Lambda =3$ is~3.

 \noindent Since the cardinal of~$V_2\setminus V_1 $ is:

$$ {\cal N}^{\cal S}_{ 2}-4= 6$$

\noindent one may note that we have the complete spectrum.

\vskip 1cm

 \subsection{Eigenvalues of~$\Delta_m$,~$m \,\in\,\N$,~$m \geq 3$}

 \noindent As previously, one can easily check that the forbidden eigenvalue~$\Lambda=2$ is not a Dirichlet one.\\

\noindent One can also check that~$ \Lambda_m =1$ and~$ \Lambda_m =3$ are eigenvalues of~$\Delta_m$.\\

 \noindent By induction, one may note that, due to the spectral decimation, the initial eigenvalue~$\Lambda_2=1$ gives birth, at this~$m^{th}$ step, to  eigenvalues~$\Lambda_{\hookrightarrow 1,m}$, and, in the same way, the initial eigenvalue~$\Lambda_2=3$ gives birth, at this~$m^{th}$ step, to eigenvalues~$\Lambda_{\hookrightarrow 3,m}$.\\
 \noindent The dimension of the Dirichlet eigenspace is equal to the cardinal of~$V_m \setminus V_{ 1}$, i.e.:

$$  {\cal N}^{\cal S}_{m } -{\cal N}^{\cal S}_{ 1 }
=3^{m }-3$$

$$
\begin{array}{|c|c|c| }
\hline \\
\text{Level}  & \text{Cardinal of the Dirichlet spectrum}     \\
\hline \\
 m   &  3^m-3 \\
\hline \\
2   &  6 \\
\hline \\
3   &   24\\
\hline \\
4   &   78 \\
\hline
\end{array}
$$

\vskip 1cm

\begin{pte}
\noindent Let us introduce:

$$\Lambda =  \displaystyle \lim_{m \to + \infty} 3^{-m}\, 4^{ m\,\delta } $$

\noindent One may note that, due to the definition of the Laplacian~$\Delta$, the limit exists.

\end{pte}

\vskip 1cm

 \subsection{Eigenvalue counting function}

 \begin{definition}\textbf{Eigenvalue counting function}\\

 \noindent Let us introduce the eigenvalue counting function, related to~$ {\cal SG}^{\cal C} \setminus V_1$, such that, for any positive number~$x$:

  $${\cal N}^{ {\cal SG}^{\cal C}\setminus V_1}(x) = \text{Card}\, \left \lbrace \Lambda \, \text{Dirichlet eigenvalue of~$-\Delta$} \, : \quad \Lambda \leq  x \right \rbrace $$

 \end{definition}

\vskip 1cm

 \begin{pte}

\noindent Given an integer~$m \geq 2$, the cardinal of~$V_m \setminus V_{ 1}$ is:

$$  {\cal N}^{\cal S}_{m } -{\cal N}^{\cal S}_{ 1 }
=3^{m }-3$$

\noindent  This leads to the existence of a strictly positive constant~$C$ such that:

$$ {\cal N}^{ {\cal SG}^{\cal C}} (C\,3^{-m}\, 4^{ m\,\delta }  )= 3^m-3 $$

 \noindent If one looks for an asymptotic growth rate of the form

$$ {\cal N}^{ {\cal SG}^{\cal C}} (x) \sim x^{\alpha_{ {\cal SG}^{\cal C}}}$$

\noindent one obtains:

$$\alpha_{ {\cal SG}^{\cal C}}=\displaystyle \frac{\ln 3}{\delta\, \ln \frac{4}{3}}=\displaystyle \frac{\ln 3}{   \frac{\ln 5}{\ln 4}\, \ln \frac{4}{3}}$$

\noindent which is not the same value as in the case of the Sierpi\'{n}ski gasket (we refer to~\cite{StrichartzLivre2006}):

$$\alpha_{ {\cal SG} }=\displaystyle \frac{\ln 3}{\ln 5}< \alpha_{ {\cal SG}^{\cal C}} $$

\noindent It appears then that increasing the number of points, and the number of connections, decreases the value of the Weyl exponent~$\alpha$.\\

\noindent By following~\cite{StrichartzLivre2006}, one may note that the ratio

$$ \displaystyle \frac{{\cal N}^{ {\cal SG}^{\cal C}} (x)}{x}$$

\noindent is bounded above and away from zero, and admits a limit along any sequence of the form~$C \,3^{-m}\, 4^{ m\,\delta }$,~\mbox{$C >0$},\\~\mbox{$m \geq 2$}. This enables one to deduce the existence of a periodic function~$g$, the period of which is equal to~$\ln \displaystyle \frac{4^{ \delta}}{3}$, discontinuous at the value~$\displaystyle \frac{4^{ \delta}}{3}$, such that:

$$\displaystyle \lim_{x \to + \infty} \left \lbrace \displaystyle \frac{{\cal N}^{ {\cal SG}^{\cal C}} (x)}{x}-g (\ln x ) \right \rbrace =0$$

 \end{pte}

\end{document}